%% file: full_paper.tex
\documentclass{IEEEtran4PSCC}

\usepackage{cite}
\usepackage{diagbox}
\usepackage{xcolor}
\usepackage{float}

\ifCLASSINFOpdf
   \usepackage[pdftex]{graphicx}
\else
   \usepackage[dvips]{graphicx}
\fi

\usepackage{multirow}
\usepackage{multicol}
\usepackage{subcaption}
\usepackage[cmex10]{amsmath}
\usepackage{amsfonts}

\usepackage{booktabs}
\usepackage{nidanfloat} 

\usepackage{url}
\usepackage{hyperref}

\makeatletter
\let\old@ps@headings\ps@headings
\let\old@ps@IEEEtitlepagestyle\ps@IEEEtitlepagestyle
\def\psccfooter#1{%
    \def\ps@headings{%
        \old@ps@headings%
        \def\@oddfoot{\strut\hfill#1\hfill\strut}%
        \def\@evenfoot{\strut\hfill#1\hfill\strut}%
    }%
    \def\ps@IEEEtitlepagestyle{%
        \old@ps@IEEEtitlepagestyle%
        \def\@oddfoot{\strut\hfill#1\hfill\strut}%
        \def\@evenfoot{\strut\hfill#1\hfill\strut}%
    }%
    \ps@headings%
}
\makeatother

\newcommand{\nonl}{\renewcommand{\nl}{\let\nl\oldnl}}

\usepackage[ruled]{algorithm2e}
\newenvironment{problem}[1][!t]
  {
   \begin{algorithm}[#1]%
  }{\end{algorithm}}

\usepackage{enumitem}

\newcommand{\PegaseSmall}{\texttt{1354\_pegase}}
\newcommand{\PegaseMedium}{\texttt{2869\_pegase}}
\newcommand{\PegaseLarge}{\texttt{9241\_pegase}}
\newcommand{\PegaseHuge}{\texttt{13659\_pegase}}
\newcommand{\SDET}{\texttt{4661\_sdet}}
\newcommand{\RTE}{\texttt{6468\_rte}}
\newcommand{\GOC}{\texttt{19402\_goc}}
\newcommand{\MSR}{\texttt{MSR}}
\newcommand{\LYON}{\texttt{LYON}}
\newcommand{\FRANCE}{\texttt{France}}

\begin{document}

\title{A Linear Outer Approximation of Line Losses for DC-based Optimal Power Flow Problems}

\author{
\IEEEauthorblockN{Haoruo Zhao, Mathieu Tanneau, Pascal Van Hentenryck}
\IEEEauthorblockA{Industrial and Systems Engineering, Georgia Institute of Technology, Atlanta, GA, United States\\
\{hzhao306, mtanneau3\}@gatech.edu, pascal.vanhentenryck@isye.gatech.edu
}
}


\maketitle

\input{tex/abstract}

\input{tex/nomenclature}

\input{tex/introduction}

\input{tex/methodology}

\input{tex/results}

\input{tex/power_flow}

\section{Conclusion} 
\label{sec:conclusion}

This work presents LLOA, a novel and simple linear model to capture line losses.
The essences of LLOA formulation are that it not only provides an intuitive and clean outer approximation to the LLQCP, which can be efficiently solved using linear programming approaches, but it also can be easily integrated within different formulations, including multi-period and stochastic formulations as well as market-clearing optimizations with discrete decisions.
Extensive numerical results comparing LLOA, LLQCP, and LLLF also demonstrate the competitive performance of LLOA against traditional methods such as LLLF, illustrating its practical usefulness.

Phase-angle formulations such as LLQCP and LLOA are less impacted by congestion than LLLF, whose performance deteriorates as the load increases.
They are also less prone to inaccuracies in line flows which, when they occur on congested lines, can lead to large increases in operating costs.
Finally, LLOA shows no drop in the ability to recover AC-feasible solutions through a power flow analysis, compared to LLQCP or LLLF.
Further numerical studies and performance benchmarks, in settings involving discrete decisions, multiple time periods or stochastic formulations, will be the topic of future work.

\section*{Acknowledgments}
This research is partly funded by NSF Award 1912244 and ARPA-E Perform Award AR0001136.

\bibliographystyle{IEEEtran}
\bibliography{IEEEabrv,references}

\newpage
\input{tex/appendix}

\end{document}

%% file: tex/abstract.tex
\begin{abstract}

This paper proposes a novel and simple linear model to capture line losses for use in linearized DC models, such as optimal power flow (DC-OPF) and security-constrained economic dispatch (SCED).
The \textit{Line Loss Outer Approximation} (LLOA) model implements an outer approximation of the line losses lazily and typically terminates in a small number of iterations.
Experiments on large-scale power systems demonstrate the accuracy and computational efficiency of LLOA and contrast it with classical line loss approaches.
The results seem to indicate that LLOA is a practical and useful model for real-world applications, providing a good tradeoff between accuracy, computational efficiency, and implementation simplicity.
In particular, the LLOA method may have significant advantages compared to the traditional loss factor formulation for multi-period, stochastic optimization problems where good reference points may not be available.
The paper also provides a comprehensive overview and evaluation of approximation methods for line losses. 

\end{abstract}
\begin{IEEEkeywords}
line losses, optimal power flow, DCOPF, SCED
\end{IEEEkeywords}

%% file: tex/nomenclature.tex
\section*{Nomenclature}

\newcommand{\Abr}{A^{\text{br}}}
\newcommand{\Agen}{A^{\text{gen}}}
\newcommand{\LF}{\text{LF}}
\newcommand{\PTDF}{\Phi}
\newcommand{\pf}{p^{\text{f}}}
\newcommand{\Pmin}{\text{P}^{\text{min}}}
\newcommand{\Pmax}{\text{P}^{\text{max}}}
\newcommand{\Tmax}{\text{T}^{\text{max}}}

\subsection*{Parameters}
    \newlist{abbrv}{itemize}{1}
    \setlist[abbrv,1]{label=,labelwidth=0.35in,align=parleft,itemsep=0.1\baselineskip,leftmargin=!}
    \begin{abbrv}
        \item[$\mathcal{N}$]    Set of buses $\mathcal{N} = \{1, ..., N\}$
        \item[$\mathcal{E}$]    Set of branches $\mathcal{E} \subseteq \mathcal{N}\times\mathcal{N}$, with $E := | \mathcal{E}|$
        \item[$\mathcal{E}^{-}$] Set of reverse branches $\mathcal{E}^{-} = \{(j, i) | (i, j) \in \mathcal{E} \}$
        \item[$\Abr$]             Branches Incidence matrix, size $E {\times} N$
        \item[$B$]                  Branch susceptance matrix, size $E {\times} E$
        \item[$r_{ij}$]                  Resistance of branch $(i, j) \in \mathcal{E}$
        \item[$\Tmax_{ij}$]            Thermal limit of branch $(i, j) \in \mathcal{E}$
        \item[$\Pmin_{i}$]            Minimum active power injection at bus $i \in \mathcal{N}$
        \item[$\Pmax_{i}$]            Maximum active power injection at bus $i \in \mathcal{N}$
        \item[$\PTDF$]               Branch-bus PTDF matrix, size $E {\times} N$
        \item[$c(\cdot)$]           Cost function for economic dispatch
    \end{abbrv}

\subsection*{Variables}
    \begin{abbrv}
        \item[$\ell^{\text{tot}}$]                   Total real power losses
        \item[$\ell_{ij}$]                 Line losses on branch $(i, j) \in \mathcal{E}$
        \item[$\overrightarrow{p_{ij}}$]        Power flow from bus $i$ to bus $j$ on branch $(i, j) \in \mathcal{E}$
        \item[$\overleftarrow{p_{ij}}$]         Power flow from bus $j$ to bus $i$ on branch $(i, j) \in \mathcal{E}$
        \item[$p_{i}$]               Active power injection at bus $i \in \mathcal{N}$
        \item[$\theta_{i}$]              Voltage angle at bus $i \in \mathcal{N}$
        \item[$v_{i}$]                   Voltage magnitude at bus $i \in \mathcal{N}$
    \end{abbrv}

%% file: tex/introduction.tex
\section{Introduction}

Optimal power flow optimization problems are ubiquitous in power system planning and operations, due to their wide applications in power blackouts prevention \cite{blackout}, power restoration \cite{Hentenryck2010VehicleRF}, and power dispatch problems \cite{Wood1984PowerGO} to name only a few.
The well-known AC optimal power flow (AC-OPF) is often regarded as the ground truth, but its nonconvexity pose significant difficulties in large-scale deployments.
The linearized DC model is often considered a reasonable approximation of AC-OPF, and is especially appealing given its formulation simplicity, computational efficiency, and robustness in many situations \cite{Stott2009DCPF}.
Security-Constrained Economic Dispatch (SCED) optimizations are widely used by independent system operators in real-time markets to produce active generator dispatches quickly.
These models typically go beyond the pure linearized DC model and capture line losses to avoid significant errors \cite{SCOPF1,SCOPF2}.

Since line losses are present during operations, it is important to study systematically how they can be captured efficiently and with high fidelity.
Many existing methods require AC-feasible prior points as warm starts \cite{linelossouter,Fitiwi2016FindingAR,linelossouter3}, which may be challenging in some settings.
Line losses can also be modeled elegantly with the linear-programming method proposed in \cite{linelossouter2}, but it was only examined on small test cases and its efficiency has not been satisfactorily evaluated. 

This paper systematically evaluates a number of approaches to model line losses in DC models.
In addition, it presents a novel and intuitive method called \textit{Line Loss Outer Approximation} (LLOA) in which quadratic line losses are iteratively approximated via linear constraints.
Iterative approximations of line losses have been discussed before (e.g., \cite{linelossouter,linelossouter2,linelossouter3}), but the LLOA method admits a simpler and clearer formulation that should be easier to deploy.
The paper evaluates the LLOA method within DC-OPF and SCED optimizations; it can also be easily integrated within many market-clearing optimizations, including those based on mixed-integer programs. The LLOA method exhibits competitive performance on many large benchmarks, both in terms of accuracy and solving times.
Being an outer approximation method, it should be ideal for multi-period, stochastic optimization problems, as well as market-clearing optimizations featuring discrete variables.
This contrasts with the traditional Loss Factor formulation where a reference point is needed to seed the method. It is also interesting to mention that the LLOA method seems to bring significant improvements for more congested test cases.

This paper is organized as follows.
Section \ref{sec:method} discusses the motivation and describes the LLOA method, including the algorithms for its application to DC-OPF and SCED.
Section \ref{sec:exp} evaluates and compares LLOA with traditional line loss approaches on a wide variety of test cases, in terms of accuracy and computational efficiency.
Section \ref{sec:power_flow} studies the recovery of AC-feasible solutions from a DC formulation.
Section \ref{sec:conclusion} concludes the paper. 

%% file: tex/methodology.tex
\section{Methods For Modeling Line Losses} 
\label{sec:method}

This section reviews the various methods to approximate line losses in linearized DC models.
These linear approximations assume that all voltage magnitudes are close to 1 and that the phase angle differences are small \cite{loss1,loss2}.
As a result, $\cos \left(\theta_{i}-\theta_{j}\right) \approx 1-\frac{\left(\theta_{i}-\theta_{j}\right)^{2}}{2}$ and the losses on branch $(i, j) \in \mathcal{E}$ can be approximated by
\begin{align}
    \label{eq:line_losses}
    \ell_{ij} & = r_{ij} \times {\overrightarrow{p_{ij}}}^{2},
\end{align}
where $r_{ij}$ is the line resistance value and $\overrightarrow{p_{ij}}$ is the active power flow on the branch.
The three methods described next differ in how Eq. \eqref{eq:line_losses} is approximated.

Throughout this section, the following conventions are used.
For simplicity, and without loss of generality, the net power injection $p_{i}$ at bus $i$ captures both generation and fixed load, if any.
The vector of active power injection at each bus is denoted by $\mathbf{p}$, whose $i^{\text{th}}$ coordinate is $p_{i}$.
Similarly, the vector of power flow (resp., reverse power flow) on each branch, with components $\overrightarrow{p_{ij}}$ (resp., $\overleftarrow{p_{ij}}$), is denoted by $\overrightarrow{\mathbf{p}}$ (resp. $\overleftarrow{\mathbf{p}}$).
Finally, $\odot$ denotes the element-wise product between two vectors.

\subsection{The DC-based Line Loss Loss Factor}
\label{sec:method:LLLF}

    The classical \textit{DC-based Line Loss Loss factor} (LLLF) approach builds on a PTDF-based DC formulation.
    Line losses are modeled via a linearization of Eq. \eqref{eq:line_losses} around a base point, which may be the latest state estimation or a previous DC solution.
    Thereby, total line losses are estimated and added to the global power balance constraint \cite{shahidehpour2003market,Litvinov2004MarginalLM,FERC}, while loss factors are introduced to distribute the losses and ensure that flows are computed accurately.
    The resulting DC model is presented in Model \ref{model:LLLF}.

    Loss factors $\LF$, loss offset $\ell^{0}$, and loss distribution factors $D$ are obtained as follows.
    First, given the reference active power injection/withdrawal $p^{\text{ref}}_{i}$ at each bus $i$, the corresponding line flows are computed using the PTDF matrix $\PTDF$, i.e.,
    \begin{align}
        {\overrightarrow{\mathbf{p}}}^{\text{ref}} = \PTDF \times \mathbf{p}^{\text{ref}}.
    \end{align}
    Then, Eq. \eqref{eq:line_losses} is differentiated with respect to line flows, and its derivative evaluated at ${\overrightarrow{\mathbf{p}}}^{\text{ref}}$, yielding the nodal loss factors
    \begin{align}
        \LF^{\top} = -2 (R \odot {\overrightarrow{\mathbf{p}}}^{\text{ref}})^{\top} \PTDF.
    \end{align}
    The loss offset $\ell^{0}$ is then defined as the scalar
    \begin{align}
        \ell^{0} = - \LF^{\top} \mathbf{p}^{\text{ref}} + \sum_{(i, j) \in \mathcal{E}} r_{ij} ({\overrightarrow{p_{ij}}}^{\text{ref}})^{2}.
    \end{align}
    The vector $D$ in the model distributes the losses among the buses, setting each element proportional to the losses on adjacent lines, i.e.,
    \begin{align}
        D_{i} = \frac{1}{2} \times \sum_{j | (i, j) \in \mathcal{E} \cup \mathcal{E}^{-}} \frac{\ell_{ij}^{\text{ref}}}{\mathbf{1}^{\top} \ell^{\text{ref}}}.
    \end{align}

    \begin{problem}
        \caption{The LLLF Model.}
        \label{model:LLLF}
        \footnotesize
        \textbf{Inputs:}
        \hspace{1in}%
        \begin{abbrv}[align=parleft,labelwidth=0.5in,leftmargin=!]
            \item[$\ \ \ell^{0}$] Loss offset
            \item[$\ \ \LF$] Loss factors
            \item[$\ \ D$] Loss Distribution factors
        \end{abbrv}
        \textbf{Variables:}
        \begin{abbrv}[align=parleft,,labelwidth=0.5in]
            \item[$\ \ \ell^{\text{tot}} \in \mathbb{R}$] Total real power losses
            \item[$\ \ \mathbf{p} \in \mathbb{R}^{N}$] Vector of active power injection at each bus
        \end{abbrv}
        \textbf{Model:}
        \begin{subequations}
        \label{eq:LLLF}
        \begin{align}
        \min \quad & c(\mathbf{p})\\
        \textrm{s.t.} \quad
            & \mathbf{1}^{\top} \mathbf{p} = \ell^{\text{tot}}\\
            & \ell^{\text{tot}} = 
                \ell^{0} + \LF^{\top}\mathbf{p}\\
            \label{eq:LLLF:thermal_limits}
            & -\Tmax \leq \PTDF \left( \mathbf{p} - \ell^{\text{tot}} D \right) \leq \Tmax\\
            & \Pmin \leq \mathbf{p} \leq \Pmax
        \end{align}
        \end{subequations}
    \end{problem}

\subsection{Line Loss Quadratic Convex Program}
\label{sec:method:LLQCP}


    \begin{problem}
        \caption{The LLQCP Model.}
        \label{model:LLQCP}
        \footnotesize
        \textbf{Inputs:}
        \hspace{1in}%
        \begin{abbrv}[align=parleft,labelwidth=0.5in,leftmargin=!]
            \item[$\ \ b_{ij}$] Susceptance of branch $(i, j) \in \mathcal{E}$
            \item[$\ \ r_{ij}$] Resistance of branch $(i, j) \in \mathcal{E}$
            \item[$\ \ s$] Slack bus index
        \end{abbrv}
        \textbf{Variables:}
        \begin{abbrv}[align=parleft,,labelwidth=0.5in]
            \item[$\ \ \theta_{i}$] Voltage phase angle at bus $i \in \mathcal{N}$
            \item[$\ \ p_{i}$] Active power injection at bus $i \in \mathcal{N}$
            \item[$\ \ \overrightarrow{p_{ij}}$] Power flow from bus $i$ to bus $j$ on branch $(i, j) \in \mathcal{E}$
            \item[$\ \ \overleftarrow{p_{ij}}$] Power flow from bus $j$ to bus $i$ on branch $(i, j) \in \mathcal{E}$
        \end{abbrv}
        \textbf{Model:}
        \begin{subequations}
        \label{eq:LLQCP}
        \begin{align}
            \label{eq:LLQCP:objective}
            \min \quad & c(\mathbf{p}) \\
            \textrm{s.t.}  \quad
            & \label{eq:LLQCP:line_flow}
                \overrightarrow{p_{ij}} = b_{ij} (\theta_{j} - \theta_{i})
                && \forall (i, j) \in \mathcal{E}\\
            & \label{eq:LLQCP:power_balance}
                p_{i} = 
                    \sum_{j | (i, j) \in \mathcal{E}^{-}} \overleftarrow{p_{ij}}
                    + \sum_{j | (i, j) \in \mathcal{E}} \overrightarrow{p_{ij}}
                && \forall{i \in \mathcal{N}}\\
            & \label{eq:LLQCP:line_losses}
                \overrightarrow{p_{ij}} + \overleftarrow{p_{ij}} \geq r_{ij} {\overrightarrow{p_{ij}}}^2  
                && \forall{(i, j) \in \mathcal{E}}\\
            & \label{eq:LLQCP:thermal_limit_from}
                {-}\Tmax_{ij} \leq \overrightarrow{p_{ij}} \leq \Tmax_{ij}
                && \forall (i, j) \in \mathcal{E}\\
            & \label{eq:LLQCP:thermal_limit_to}
                {-}\Tmax_{ij} \leq \overleftarrow{p_{ij}} \leq \Tmax_{ij}
                && \forall (i, j) \in \mathcal{E}\\
            & \label{eq:LLQCP:generation_limit}
                \phantom{-}\Pmin_{i} \leq \mathbf{p} \leq \Pmax_{i}
                && \forall i \in \mathcal{N}\\
            & \label{eq:LLQCP:slack_bus}
                \theta_{s} = 0
        \end{align}
        \end{subequations}
    \end{problem}

\begin{figure}
        \centering
        \includegraphics[width=\columnwidth]{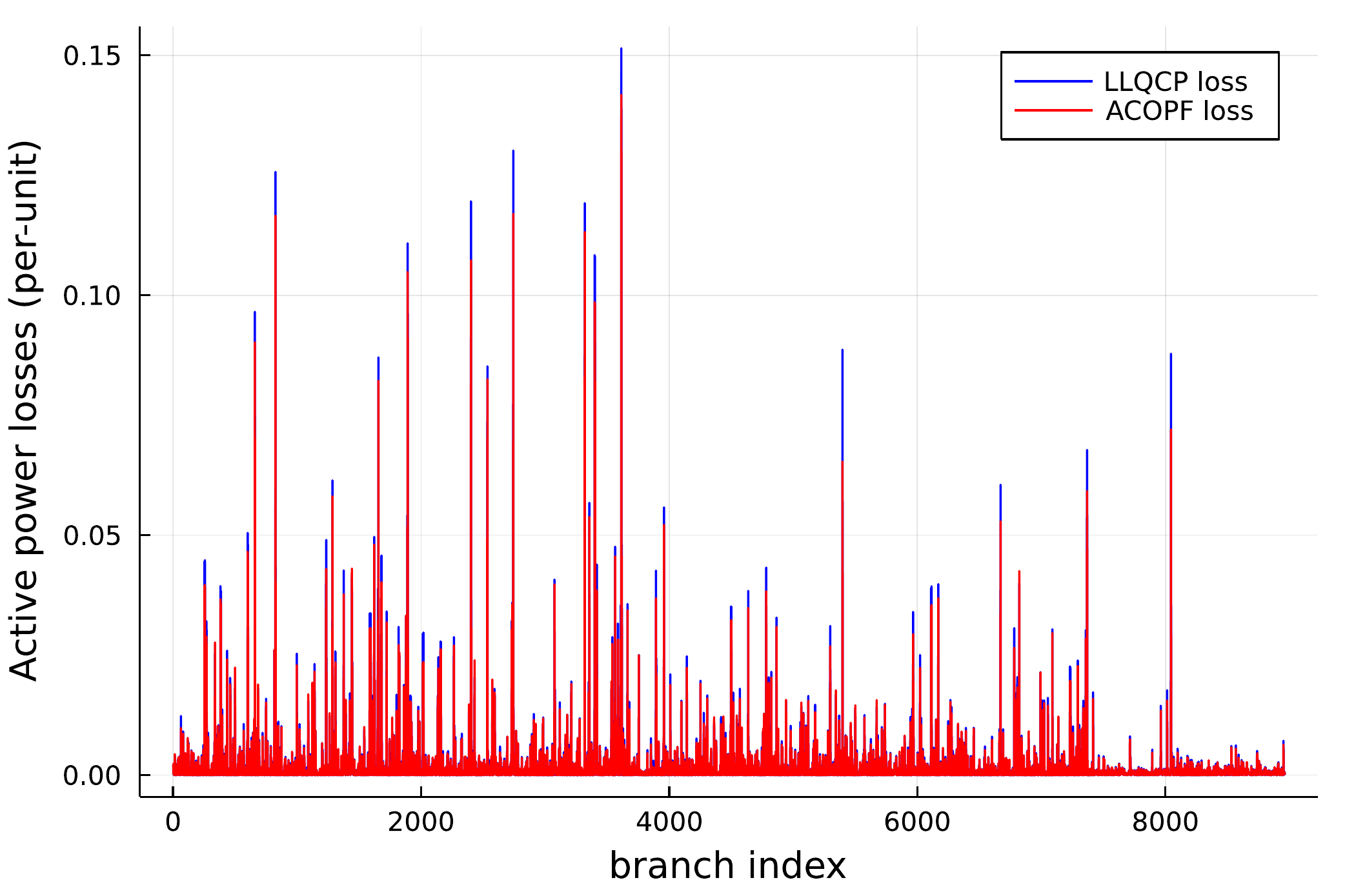}
        \caption{Real vs. Estimated DCOPF Line Losses on the France RTE system using LLQCP.}
        \label{fig:LLQCP:real_vs_est_loss}
    \end{figure}

    The \emph{Line Loss Quadratic Convex Program} (LLQCP) \cite{Coffrin2012ApproximatingLL} is a  DC formulation based on phase angles.
    It captures the flow direction on each branch $(i, j) \in \mathcal{E}$ through two variables $\overrightarrow{p_{ij}}$ and $\overleftarrow{p_{ij}}$ and adds the following quadratic constraint on each branch:
    \begin{align}
        \overrightarrow{p_{ij}} + \overleftarrow{p_{ij}} \geq \ell_{ij} = r_{ij} (\overrightarrow{p_{ij}})^2.
    \end{align}
    The formulation is specified in Model \ref{model:LLQCP}.
    Constraint \eqref{eq:LLQCP:power_balance} is the nodal power balance at bus $i \in\mathcal{N}$, which states an equivalence between nodal net active power injection $p_{i}$ and the sum of active power flowing in ($\overleftarrow{p_{ij}}$) and out ($\overrightarrow{p_{ij}}$) of the bus.
    Line losses are modeled via constraint \eqref{eq:LLQCP:line_losses}, and thermal limits are enforced by constraints \eqref{eq:LLQCP:thermal_limit_from}-\eqref{eq:LLQCP:thermal_limit_to}.
    Other constraints are identical to the vanilla DC formulation.
    Figure \ref{fig:LLQCP:real_vs_est_loss} compares the real and estimated line losses on the French system and it highlights the fact that the LLQCP model captures the line losses with high fidelity.
     
\begin{figure}
        \centering
        \includegraphics[width=\columnwidth]{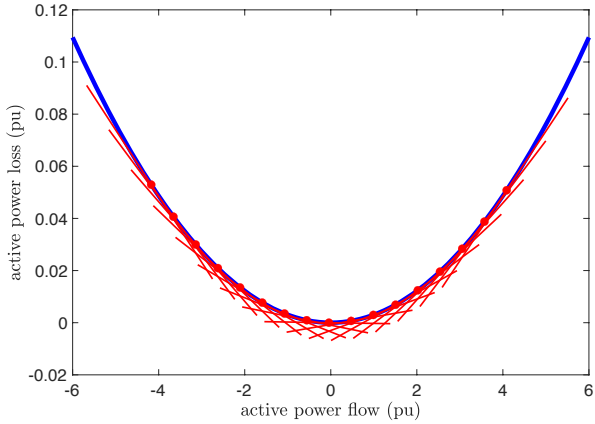}
        \caption{The LLQCP Linearization considered in \cite{Coffrin2012ApproximatingLL}. The quadratic relation of Eq. \eqref{eq:line_losses} is depicted in \textcolor{blue}{\textbf{blue}}, and the static linearization in \textcolor{red}{\textbf{red}}.}
        \label{fig:LLQCP}
    \end{figure}

    The main drawback of the LLQCP formulation is the nonlinearity of constraint \eqref{eq:LLQCP:line_losses}, which leads to increased computing times and can make the model challenging when discrete decisions are included.
    Thus, in the original presentation \cite{Coffrin2012ApproximatingLL}, the authors consider a static linearization as illustrated in Figure \ref{fig:LLQCP}.
    Nevertheless, computational and accuracy issues were reported on large test cases.
     
\subsection{The Line Loss Outer Approximation}
\label{sec:method:LLOA}

    To reduce the computational burden and maintain the accuracy of the LLQCP model, the LLOA approach replaces the  static linearization of the quadratic terms with a lazy outer approximation guided by a candidate solution to the model.
    Given a good approximation  $\overrightarrow{p_{ij}}^{\text{ref}}$ of the active power flow on branch $(i, j)$, constraint \eqref{eq:LLQCP:line_losses} can be outer-approximated with the \emph{linear} constraint 
    \begin{align}
        \label{eq:LLOA:line_losses}
        \overrightarrow{p_{ij}} + \overleftarrow{p_{ij}}
        \geq 
        r_{ij} \left[
            -(\overrightarrow{p_{ij}}^{\text{ref}})^2 + 2 \times \overrightarrow{p_{ij}}^{\text{ref}} \times \overrightarrow{p_{ij}}
        \right].
    \end{align}
    More generally, given a set of $H$ linearization points $\mathcal{F} = \{\overrightarrow{\mathbf{p}}^{1}, ..., \overrightarrow{\mathbf{p}}^{H}\}$, the LLOA formulation includes, for each branch, $H$ constraints of the form \eqref{eq:LLOA:line_losses}.
    The corresponding formulation is given in Model \ref{model:LLOA}: it is identical to Model \ref{model:LLQCP}, except for constraint \eqref{eq:LLQCP:line_losses} which is replaced with outer-approximation constraints of the form \eqref{eq:LLOA:line_losses}.


    \begin{problem}
        \caption{The LLOA Model.}
        \label{model:LLOA}
        \footnotesize
        \textbf{Inputs:}
        \begin{abbrv}[align=parleft,labelwidth=0.5in,leftmargin=!]
            \item[$\ \ b_{ij}$] Susceptance of branch $(i, j) \in \mathcal{E}$
            \item[$\ \ r_{ij}$] Resistance of branch $(i, j) \in \mathcal{E}$
            \item[$\ \ s$] Slack bus index
            \item[$\ \ \mathcal{F}$] Set of linearization points
        \end{abbrv}
        \textbf{Variables:}
        \begin{abbrv}[align=parleft,,labelwidth=0.5in]
            \item[$\ \ \theta_{i}$] Voltage phase angle at bus $i \in \mathcal{N}$
            \item[$\ \ p_{i}$] Active power injection at bus $i \in \mathcal{N}$
            \item[$\ \ \overrightarrow{p_{ij}}$] Power flow from bus $i$ to bus $j$ on branch $(i, j) \in \mathcal{E}$
            \item[$\ \ \overleftarrow{p_{ij}}$] Power flow from bus $j$ to bus $i$ on branch $(i, j) \in \mathcal{E}$
        \end{abbrv}
        \textbf{Model:}
        \begin{align*} 
            \min \quad & c(\mathbf{p}) \\
            \text{s.t.} \quad 
            & \eqref{eq:LLQCP:line_flow} \text{ -- } \eqref{eq:LLQCP:power_balance}\\
            & \eqref{eq:LLOA:line_losses} 
                && \forall (i, j) \in \mathcal{E}, \forall \overrightarrow{\mathbf{p}}^{\text{ref}} \in \mathcal{F}\\
            & \eqref{eq:LLQCP:thermal_limit_from} \text{ -- } \eqref{eq:LLQCP:slack_bus}
        \end{align*}
    \end{problem}

    The general LLOA iterative algorithm is described in Algorithm \ref{alg:LLOA}, in the context of DC-OPF.
    Besides the power network data, which is omitted for brevity, the only input is a convergence tolerance $\epsilon > 0$: unlike LLLF, no reference solution is needed for correctness.
    At each iteration, the current approximation is solved and, if the convergence criterion is not met, the outer approximation is refined using the obtained solution.
    This process is equivalent to applying Kelley's cutting-plane algorithm to the LLQCP formulation, thereby yielding, in the limit, the same accuracy as solving LLQCP.
    Furthermore, the LLOA model remains linear throughout the procedure, it has a smaller size compared to the static linearization used in \cite{Coffrin2012ApproximatingLL}, and the algorithm can be terminated at any point to produce an appropriate trade-off between efficiency and accuracy.
    Finally, note that, although a starting point is not required for correctness, the algorithm can be warm-started using a reference solution, e.g., a previous operating point, or the solution of the vanilla DCOPF.


    \begin{algorithm}[!t]
        \SetAlgoLined
        \caption{The LLOA iterative algorithm.} \nonl
        \label{alg:LLOA}
        \footnotesize
        \textbf{Inputs:}
        \begin{abbrv}[align=parleft,labelwidth=0.75in,leftmargin=!]
            \item[$\ \ \epsilon$]   Convergence tolerance
        \end{abbrv}
        \textbf{Algorithm:}\\
        Set $\mathcal{F} = \emptyset, \Delta^{0} = +\infty, Z^{0} = +\infty, k=0$\\
        \While{$|\Delta^{k}| > \epsilon$}{
            Solve current approximation using Model \ref{model:LLOA}, and retrieve\\
            \hspace{0.1in} corresponding dispatch $\mathbf{p}^{k}$ and flow $\overrightarrow{\mathbf{p}}^{k}$.\\
            $Z^{k} = c\left( \mathbf{p}^{k} \right)$\\
            $\Delta^{k} = \dfrac{|Z^{k} - Z^{k-1}|}{|Z^{k-1}|}$\\
            $\mathcal{F} \longleftarrow \mathcal{F} \cup \{ \overrightarrow{\mathbf{p}}^{k} \}$\\
            $k \longleftarrow k + 1$\\
        
        }%
        \textbf{Return:}
        \begin{abbrv}[align=parleft,labelwidth=0.75in,leftmargin=!]
            \item[$\ \ \mathbf{p}^{*}$]   Optimal active power dispatch
            \item[$\ \ Z^{*} = c(\mathbf{p}^{*})$]   Optimal objective value
        \end{abbrv}
    \end{algorithm}

\subsection{Discussion}
\label{sec:method:discussion}

While the three methods LLLF, LLQCP, and LLOA all adjust the vanilla DC formulation to approximate line losses based on Eq. \eqref{eq:line_losses}, they exhibit some fundamental differences, whose benefits and limitations are discussed next.

First, on the one hand, the LLLF method is a top-down approach: total losses are estimated using a reference point, then distributed onto the lines and buses in order to account for thermal limits and power balance.
In particular, the LLLF method only differs from a vanilla PTDF-based DC formulation by the addition the total line loss variable $\ell^{\text{tot}}$, and the modification of PTDF coefficients using loss distribution factors.
This approach is best-suited to low-congestion settings where few, albeit dense, PTDF constraints are needed.
On the other hand, LLQCP and LLOA are bottom-up approaches, wherein losses are estimated on each branch through the convexification of Eq. \eqref{eq:line_losses}, without the need for a reference point.
These methods rely on sparser phase-angle formulations, at the price of a higher number of variables.

Second, the LLOA formulation can be dynamically refined to improve its accuracy as needed, which exploits the warm-starting capabilities of simplex-based LP solvers.
In contrast, because the loss factor coefficients are embedded in the problem's constraint matrix, the LLLF formulation cannot be refined without re-building the problem, effectively forcing a cold-start subsequent solve.
While doing so may be tractable for single-period, deterministic problems, this is no longer the case in settings where re-optimizing a model from scratch is costly, for instance, when considering multi-period and/or stochastic problems, or when discrete decisions are involved, e.g., in unit-commitment problems.

Third, in contrast to LLLF, LLQCP and LLOA do not require an initial solution for correctness.
While a reasonably good reference point is typically available for solving real-time, single-period OPF problem, this is no longer the case when considering multiple time periods spanning several hours in the future, e.g., in day-ahead markets.
This limitation is made worse by the fact that longer-horizon problems are also more expensive to re-solve if the loss factors in LLLF were to be updated.

Finally, while presented in the context of OPF problems, all methods can obviously be applied to the Security-Constrained Economic Dispatch (SCED) \cite{MISO2009_SCED,BPM_002_D}; note that the SCED formulation used by MISO in its real-time market is based on LLLF.
The inclusion of LLQCP or LLOA formulations only impacts the representation of line flows, namely, using a phase angle-based formulation instead of the PTDF approach used by MISO.
Once this formulation is available, the line-loss constraints are similar to those of the DC model just presented, and the rest of the model is unchanged.

%% file: tex/results.tex
\section{Computational Experiments}
\label{sec:exp}

This section compares the above formulations on a collection of large-scale test cases.
All models are formulated in Julia 1.6 using JuMP \cite{DunningHuchetteLubin2017}; all DC-based formulations are solved using Gurobi 9.1.0 \cite{gurobi} with default parameters, and all AC models are solved with Ipopt.
Experiments are carried out on Linux machines with dual Intel Xeon 6226@2.7GHz CPUs on the PACE Phoenix cluster \cite{PACE}.
Section \ref{sec:res:testcases} presents the test cases.
Section \ref{sec:res:LLOA_vs_LLQCP} studies the trade-offs between LLOA and LLQCP, and Section \ref{sec:res:accuracy} compares the accuracy of the considered DC-based approaches against an AC formulation.
Computing times are reported in Section \ref{sec:res:computing_times}

\subsection{Test cases}
\label{sec:res:testcases}

    Experiments are carried out on a collection of large-scale instances from PGLib \cite{babaeinejadsarookolaee2021power}, as well as three instances provided by the French transmission system operator RTE.
    The size of the test cases are given in Table \ref{tab:res:systems}, which reports the number of buses, branches, and generators in each system.
    The three test cases from RTE are the following: $\MSR$ is a snapshot of a south-east region in France known to suffer from voltage issues, $\LYON$ is the French high-voltage system plus low- and medium-voltage grid components in the Lyon area, and $\FRANCE$ is the full French transmission grid.

    For each system, both a vanilla OPF and MISO's real-time SCED formulations are considered; the latter is described in \cite{MISO2009_SCED,BPM_002_D}.
    The SCED formulation mainly differs from OPF in the inclusion of reserve requirements and contingency constraints.
    In the current formulation used by MISO, line losses are modeled via loss factors, as in LLLF, and reserve requirements, power balance, and transmission constraints are soft, i.e., they can be violated at a (high) penalty.
    The AC-SCED formulation considers the same modeling of reserve requirements, but uses explicit AC power flow equations.
    Consequently, in order to keep the AC-SCED formulation tractable, only a single generator contingency is considered, corresponding to the loss of the largest generator.
    Finally, all models are expressed in per-unit, with a base scaling factor of 100MW (resp. 100MVAR) for active (resp. reactive) power.
    
    All DC-based formulations are compared to the AC formulation using the following metrics: the percent objective gap with respect to the AC solution found by Ipopt, the total estimated line losses, and the Mean Absolute Error (MAE) of active power dispatch.
    Note that, because DC-based formulations are \emph{approximations}, the objective gap may be positive (higher than AC) or negative (lower than AC).
    Letting $\mathbf{pg}^{*}$ and $\mathbf{pg}_{AC}^{*}$ denote the vectors of active power dispatch for the considered DC-based and AC formulation, respectively, the MAE is given by
    \begin{align*}
        \text{MAE} = \frac{1}{G} \left\| \mathbf{pg}^{*} - \mathbf{pg}^{*}_{AC} \right\|_{1},
    \end{align*}
    where $G$ is the number of generators.

    \begin{table}[!t]
        \centering
        \begin{tabular}{llrrr}
            \toprule
            Source & System & Buses &  Branches & Generators \\
            \midrule
            PGLib 
            & \PegaseSmall  &   1354 &  1991 &  260 \\
            & \PegaseMedium &   2869 &  4582 &  510 \\
            & \SDET         &   4661 &  5997 & 1176 \\
            & \RTE          &   6468 &  9000 & 1295 \\
            & \PegaseLarge  &   9241 & 16049 & 1445 \\
            & \PegaseHuge   &  13659 & 20467 & 4092 \\
            & \GOC          &  19402 & 34704 &  971 \\
            \midrule
            RTE 
            & \MSR          &    403 &   550 &  115 \\
            & \LYON         &   3411 &  4499 &  771 \\
            & \FRANCE       &   6705 &  8962 & 1708 \\
            \bottomrule
        \end{tabular}
        \caption{Test Cases from PGLib and RTE}
        \label{tab:res:systems}
    \end{table}

\subsection{Tradeoff between LLOA and LLQCP}
\label{sec:res:LLOA_vs_LLQCP}

\begin{figure}[!t]
        \centering
        \includegraphics[width=\linewidth]{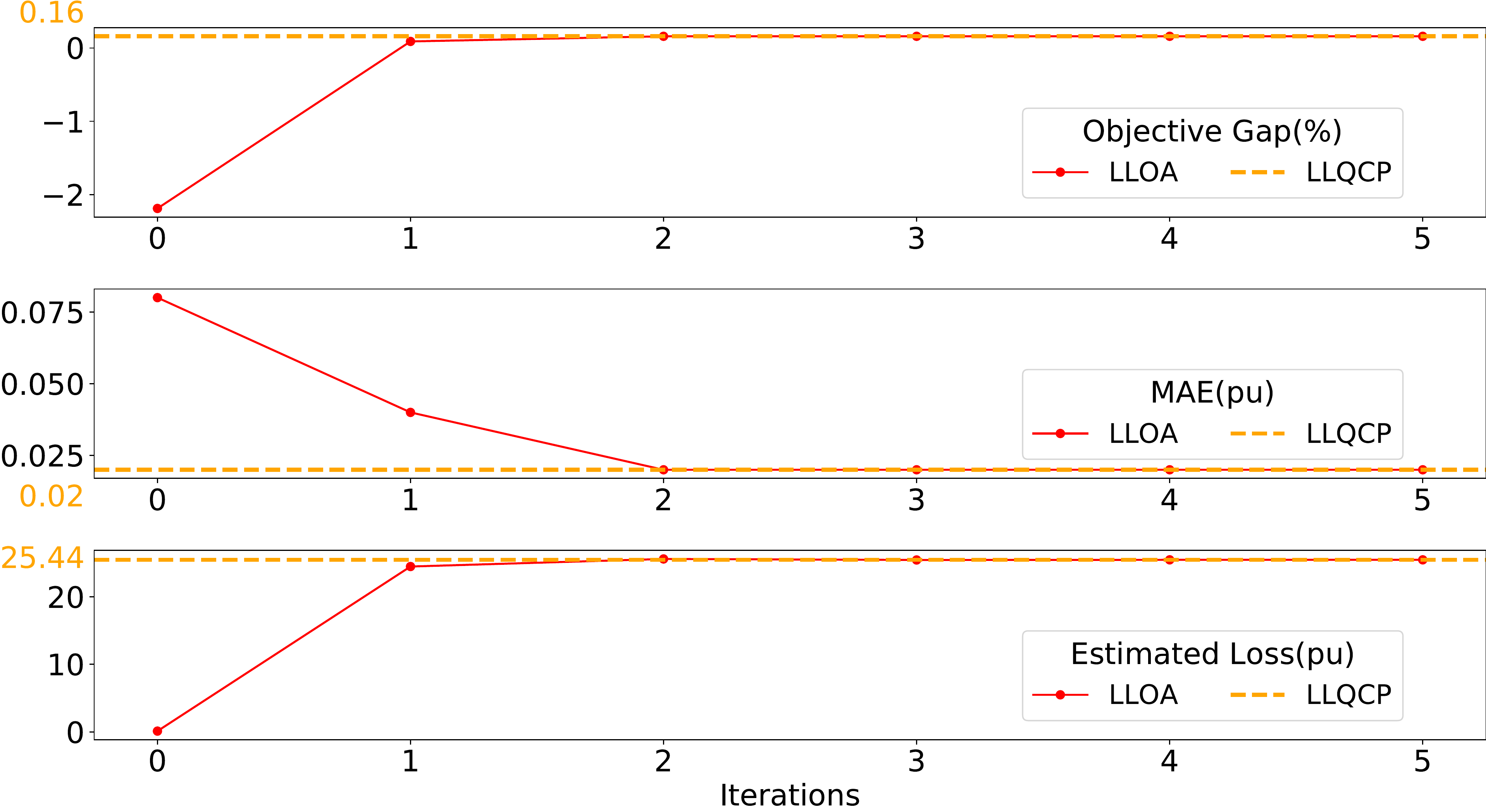}
        \caption{Convergence of LLOA on the \PegaseMedium\ system.}
        \label{fig:res:LLOA:convergence}
    \end{figure}

    The convergence of Kelley's cutting-plane algorithm in LLOA is illustrated in Figure \ref{fig:res:LLOA:convergence} for the  \PegaseMedium\ system, with a strict convergence tolerance $\epsilon = 10^{-6}$.
    Values obtained for the LLQCP solution are also displayed as a reference.
    The behavior observed in Figure \ref{fig:res:LLOA:convergence} is representative of LLOA's behavior across the overall testset.
    Indeed, across all considered instances, LLOA always reaches the default $\epsilon = 10^{-3}$ tolerance in no more than 3 iterations.
    As can be seen in Figure \ref{fig:res:LLOA:convergence}, further iterations do not improve the approximation significantly.
    
    Note that, if no initial reference point is provided, then the first iteration (iteration $0$ in Figure \ref{fig:res:LLOA:convergence}) is equivalent to solving a vanilla DC formulation, with no line losses.
    To enable a fair comparison to LLLF, which requires a reference solution, in all that follows, LLOA is seeded with the same initial solution as LLLF, obtained by solving the DC formulation with no losses.
    Accordingly, the number of iterations required by LLOA to converge is typically 1 and never more than 2.
    Therefore, subsequent experiments use a single LLOA iteration, and the corresponding method is denoted by LLOA$_{1}$.
    
    
\subsection{Accuracy Comparison}
\label{sec:res:accuracy}

    Tables \ref{tab:res:OPF:objective_gap}-\ref{tab:res:OPF:losses} report, for each method and test system, the percent objective gap, the MAE of active power dispatch, and the total estimated line losses, for the OPF problems.
    The same statistics are reported for SCED problems in Tables \ref{tab:res:SCED:objective_gap}-\ref{tab:res:SCED:losses}.
    In each table, the best-performing method is indicated in bold.
    Because the vanilla DC formulation (denoted by DC in the tables) does not model line losses, it is removed from Tables \ref{tab:res:OPF:losses} and \ref{tab:res:SCED:losses} which, instead, report total line losses in the AC solution.
    The results are consistent across OPF and SCED problems, and motivate the following comments.
    
    \begin{table}[!t]
        \centering
        \begin{tabular}{lrrrr}
            \toprule
            \textbf{Obj Gap(\%)} & DC & LLOA$_1$   & LLQCP & LLLF   \\
            \hline
            \PegaseSmall &  -3.23 &   +0.26 &   +0.30 &   \textbf{+0.09} \\
            \PegaseMedium &  -3.11 &  -0.32 &   \textbf{+0.24} &  -0.50 \\
            \SDET    & -1.56 & +0.05  & +0.09  & \textbf{-0.04} \\
            \RTE     & -4.20 & -0.46 & \textbf{-0.45} & -0.82 \\
            \PegaseLarge  & -3.46 & \textbf{-0.21}  & +0.21  & -0.42 \\
            \PegaseHuge & -2.07 & \textbf{-0.07}  & +0.08  & +0.13  \\
            \GOC    & -4.04 & \textbf{-0.02} & +0.03  & +0.06  \\
            \MSR          & -0.14 & +0.01  & +0.01  & \textbf{+0.01}  \\
            \LYON         & -1.79 & \textbf{+0.85}  & +0.90  & +1.04  \\
            \FRANCE       & -3.64 & -0.18 & \textbf{-0.15} & -0.21 \\
            \bottomrule
        \end{tabular}%
        \caption{Percent Objective Gap (OPF).}
        \label{tab:res:OPF:objective_gap}
    \end{table}

    \begin{table}[!t]
        \centering
        \begin{tabular}{lrrrr}
            \toprule
            \textbf{MAE(pu)} & DC & LLOA$_1$   & LLQCP & LLLF   \\
            \hline
            \PegaseSmall  & 0.08 & 0.02 & \textbf{0.02} & 0.04 \\
            \PegaseMedium  & 0.13 & 0.19 & \textbf{0.02} & 0.16 \\
            \SDET    & 0.04 & 0.04 & \textbf{0.03} & 0.04 \\
            \RTE     & 0.10 & 0.05 & \textbf{0.05} & 0.06 \\
            \PegaseLarge  & 0.12 & 0.14 & \textbf{0.02} & 0.13 \\
            \PegaseHuge & 0.04 & 0.03 & \textbf{0.01} & 0.04 \\
            \GOC    & 0.04 & 0.00 & \textbf{0.00} & 0.00 \\
            \MSR           & 0.00 & 0.00 & \textbf{0.00} & 0.00 \\
            \LYON          & 0.04 & \textbf{0.01} & 0.01 & 0.01 \\
            \FRANCE        & 0.01 & 0.01 & \textbf{0.00} & 0.01 \\
            \bottomrule
        \end{tabular}
       \caption{Active Power Dispatch MAE (OPF).}
        \label{tab:res:OPF:MAE}
    \end{table}
    
    \begin{table}[!t]
        \centering
        \begin{tabular}{lrrrr}
            \toprule
            \textbf{Estimated Loss(pu)} & LLOA$_1$   & LLQCP & LLLF & AC  \\
            \hline
            \PegaseSmall  & 15.99  & 16.19 & \textbf{15.46} & 15.03 \\
            \PegaseMedium  & 19.82  & \textbf{29.67} & 18.61 & 27.02 \\
            \SDET    & 11.70  & \textbf{12.29} & 11.34 & 12.20 \\
            \RTE     & 20.42  & 20.55 & \textbf{18.60} & 19.45 \\
            \PegaseLarge  & 55.87  & \textbf{77.06} & 53.88 & 72.26 \\
            \PegaseHuge & 73.88 & \textbf{90.73} & 80.93 & 88.71 \\
            \GOC    & 35.68  &\textbf{36.29} & 36.34 & 36.18 \\
            \MSR           & 0.50   & 0.50  & \textbf{0.50}  & 0.47  \\
            \LYON          & 12.73  & 12.89 & \textbf{10.50} & 10.10 \\
            \FRANCE        & \textbf{12.77}  & 13.67 & 12.60 & 13.14 \\
            \bottomrule
        \end{tabular}
        \caption{Total estimated line losses (OPF).}
        \label{tab:res:OPF:losses}
    \end{table}

    \begin{table}[!t]
        \centering
        \begin{tabular}{lrrrr}
            \toprule
            \textbf{Obj Gap(\%)} &   DC &   LLOA$_1$ &  LLQCP &   LLLF \\
            \midrule
            \PegaseSmall &  -2.31 &   +0.28 &   +0.29 &   \textbf{+0.15} \\
            \PegaseMedium &  -2.18 &   \textbf{+0.02} &   +0.16 &  -0.13 \\
            \SDET &   -1.50 &   +0.07 &   +0.09 & \textbf{ -0.01 }\\
            \RTE &  -3.77 &  -0.16 & \textbf{ -0.12} &  -0.49 \\
            \PegaseLarge &  -3.05 &  -0.18 &   \textbf{+0.12} &  -0.26 \\
            \PegaseHuge &  -2.03 &  -0.14 &  \textbf{ +0.06} &   +0.07 \\
            \GOC &  -3.91 &  -0.09 &  -0.04 &  \textbf{-0.02} \\
            \MSR &  -4.58 &   +0.21 &   +0.25 &   \textbf{+0.16} \\
            \LYON &   -1.80 &   \textbf{+0.81} &   +0.83 &   +1.05 \\
            \FRANCE &  -3.55 &   -0.19 &  \textbf{-0.15} &  -0.24 \\
            \bottomrule
        \end{tabular}
        \caption{Percent Objective Gap (SCED).}
        \label{tab:res:SCED:objective_gap}
    \end{table}

    \begin{table}[!t]
        \centering
        \begin{tabular}{lllll}
        \toprule
               \textbf{MAE(pu)} &  DC &  LLOA$_1$ & LLQCP &  LLLF \\
                \midrule
                \PegaseSmall &  0.07 &  0.03 &  0.03 &  \textbf{0.02} \\
                \PegaseMedium &  0.09 &  0.07 &  \textbf{0.02} &  0.06 \\
                \SDET &  0.03 &  0.03 &  \textbf{0.03} &  0.03 \\
                \RTE &  0.12 &  0.08 &  \textbf{0.07} &  0.08 \\
                \PegaseLarge &  0.13 &  0.11 &  \textbf{0.05} &  0.11 \\
                \PegaseHuge &  0.04 &  0.04 &  \textbf{0.01} &  0.05 \\
                \GOC &  0.04 &  0.01 &  \textbf{0.01} &  0.01 \\
                \MSR &  0.01 &   0.00 &   \textbf{0.00} &   0.00 \\
                \LYON &  0.02 &  0.02 &  \textbf{0.01} &  0.03 \\
                \FRANCE &  0.01 &  0.01 &   \textbf{0.00} &  0.01 \\
        \bottomrule
        \end{tabular}
        \caption{Active Power Dispatch MAE (SCED).}
        \label{tab:res:SCED:MAE}
    \end{table}

    \begin{table}[!t]
        \centering
        \begin{tabular}{lrrrr}
            \toprule
            \textbf{Estimated Loss(pu)} &   LLOA$_1$ &  LLQCP &   LLLF &  AC \\
            \midrule
            \PegaseSmall &  15.47 &  15.53 &  \textbf{14.78} &   14.51 \\
            \PegaseMedium &  \textbf{22.77} &  25.44 &  21.04 &   23.69 \\
            \SDET &  11.65 &  \textbf{12.06} &  11.06 &   11.98 \\
            \RTE &  \textbf{19.43} &  19.89 &  17.66 &   18.92 \\
            \PegaseLarge &   55.90 &  \textbf{74.39} &  54.14 &  71.54 \\
            \PegaseHuge &  68.11 &  \textbf{89.58} &  74.39 &  87.77 \\
            \GOC &  36.66 &  \textbf{37.41} &  37.28 &  37.38 \\
            \MSR &   0.46 &   0.46 &   \textbf{0.45} &   0.44 \\
            \LYON &  13.79 &   13.90 &   \textbf{10.10} &  10.46 \\
            \FRANCE &  12.29 &  \textbf{12.76} &  12.15 &  12.59 \\
            \bottomrule
        \end{tabular}%
        \caption{Total estimated line losses (SCED).}
        \label{tab:res:SCED:losses}
    \end{table}
    
    

    First, LLLF, LLQCP, and LLOA\textsubscript{1} all outperform the vanilla DC formulation in terms of objective gap, MAE, and total estimated losses.
    Indeed, because it ignores line losses, the vanilla DC formulation is found to systematically under-estimate the total objective cost, by as much as $4.20\%$ on OPF and $4.58\%$ on SCED problems.
    Furthermore, ignoring line losses results in lower total active power generation which, in unit-commitment problems, may lead to committing too few units, thereby causing price spikes and/or reliability issues later in the day.
    Since all electricity markets in the US are based on DC formulations, this observation reinforces the fact that line loss information should be included in TSO market specifications.
    
    Second, among LLLF, LLQCP, and LLOA\textsubscript{1}, no individual method dominates the others across all considered instances.
    Nevertheless, the LLQCP method is found to be the most accurate overall, with typically smaller objective gaps and MAE, as well as better line loss estimation.
    This demonstrates the benefits of exploiting nonlinear information in the formulation although, as will be shown next, this comes at the price of higher computing times.
    In addition, it establishes the accuracy that LLOA can reach since, in the limit, it solves the same formulation as LLQCP.
    
\begin{figure*}[!t]
        \centering
        \includegraphics[width=0.8 \linewidth]{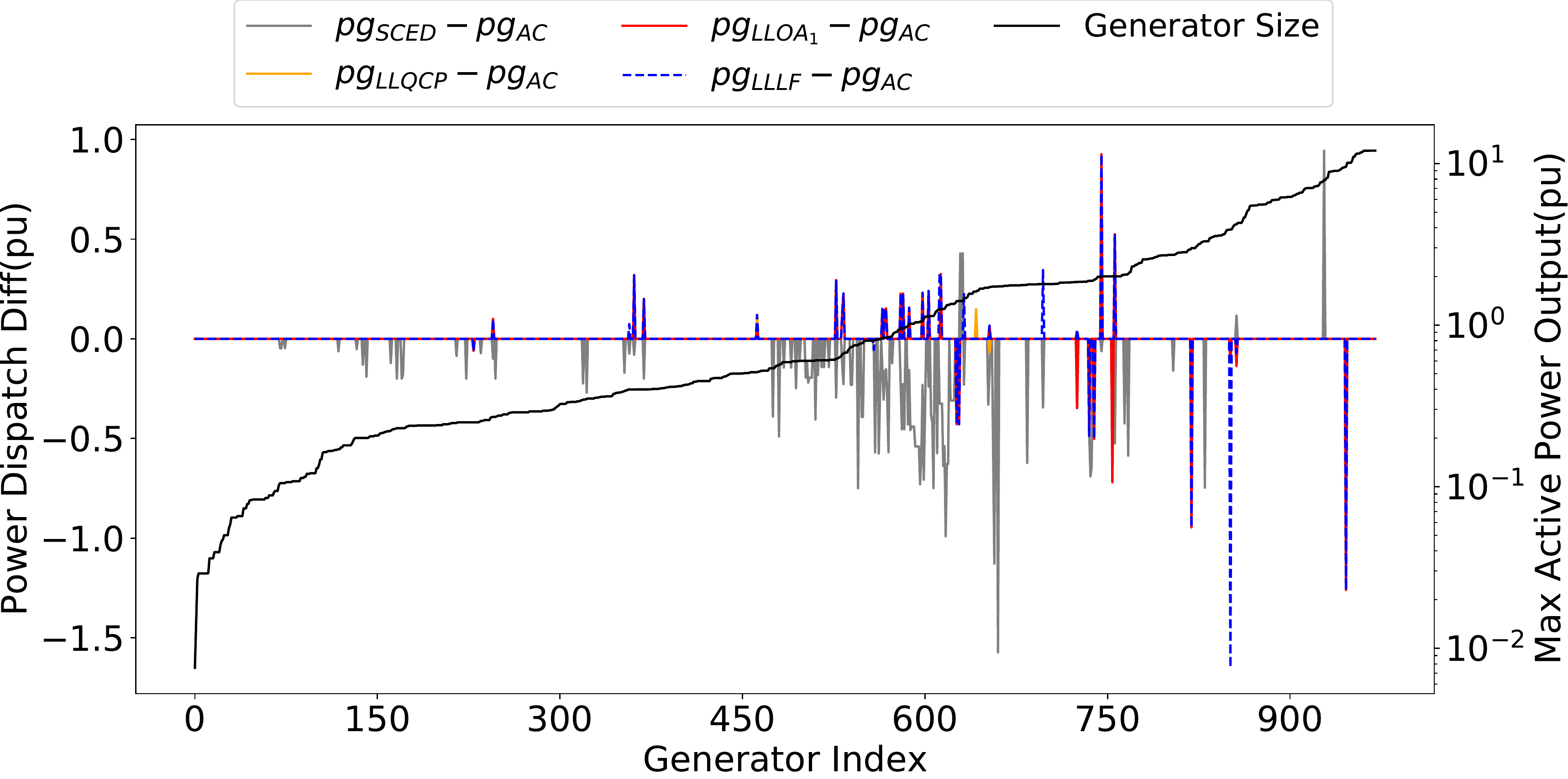}
        \caption{%
        Differences in active power dispatch (left axis, in p.u.) between DC-based and AC SCED formulations, on the \GOC\ system. Generators are ordered from left to right by increasing maximum capacity (right axis, in p.u.).
        }
        \label{fig:rel_gen}
    \end{figure*}

\begin{figure*}[!t]
        \centering
        \includegraphics[width=\linewidth]{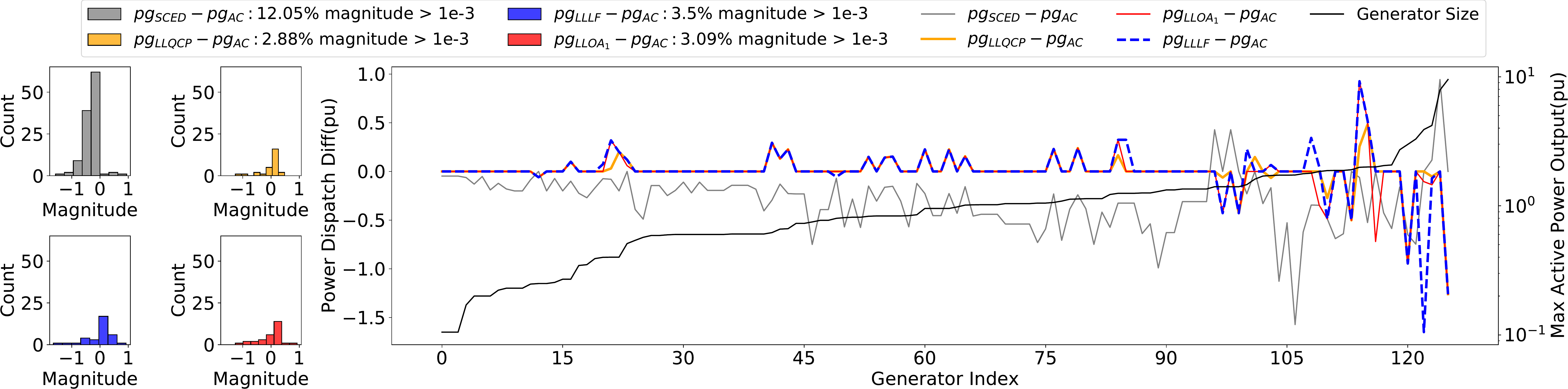}
        \caption{%
        Histogram (left) and individual differences (right) of active power dispatch (left axis, in p.u.) between DC-based and AC SCED formulations, on the \GOC\ system.
        Only generators for which at least one method yields an error ${>}10^{-3}$p.u. are displayed, in increasing order of maximum capacity (right axis, in p.u.).
        }
        \label{fig:rel_gen_threshold}
    \end{figure*}

    In light of the above observations, it is interesting to look at each method's behavior at the individual generator level.
    To that end, Figure \ref{fig:rel_gen} plots the difference in active power dispatch between each DC-based method, and the reference AC solution, for the SCED problem on the \GOC\ system.
    Individual differences are depicted (left axis), together with the maximum capacity of each generator (right axis, in log-scale).
    Figure \ref{fig:rel_gen_threshold} presents the same data after filtering out generators for with small (${<}0.1$MW) differences, as well as the corresponding histograms of active power dispatch differences, for each DC-based method.
    Unsurprisingly, the generator dispatch in the vanilla DC solution are usually lower than in the AC solution.
    Moreover, with about $12\%$ of generators with non-zero difference, the vanilla DC solution differs from the AC solution significantly more than LLLF, LLQCP, and LLOA\textsubscript{1}, which display both fewer and smaller differences.
    Remarkably, the LLLF, LLQCP, and LLOA\textsubscript{1} solutions are very similar, with about $97\%$ of generators having an almost identical dispatch to the AC solution.
    This latter behavior is highly desirable, since it means that operators need only adjust a handful of generators to recover an AC-feasible operating point; this question is further investigated in Section \ref{sec:power_flow}.
    

\subsection{Computational Efficiency}
\label{sec:res:computing_times}

    Tables \ref{tab:res:OPF:computing_time} and \ref{tab:res:SCED:computing_time} report the computing time of each method on OPF and SCED problems.
    The reported times are in seconds, and do not include model setup time, which is negligible, nor the computation of PTDF coefficients, which can be performed offline as it only depends on the network topology.
    Finally, computing times for LLLF and LLOA\textsubscript{1} do not include the time needed to obtain the initial reference solution since, in real-time operations, the previous solution or current state estimation is readily available.


    \begin{table}[!t]
        \centering
        \begin{tabular}{lrrrrr}
            \toprule
            \textbf{Solving time (secs)} & DC & LLOA$_1$ &  LLQCP &  LLLF &  AC\\
            \midrule
                \PegaseSmall    &  0.03 &   0.31 &   0.29 &  0.03 &  25.71 \\
                \PegaseMedium   &  0.14 &   1.01 &   0.83 &  0.05 &  41.48 \\
                \SDET           &  1.75 &   2.06 &   1.04 &  1.04 &  47.28 \\
                \RTE            &  0.02 &   3.41 &   2.75 &  0.02 &  161.91 \\
                \PegaseLarge    &  0.60 &   5.10 &   5.42 &  0.34 &  106.56 \\
                \PegaseHuge     &  2.24 &   6.83 &   6.74 &  0.99 &  141.20 \\
                \GOC            &  0.04 &   4.34 &  20.49 &  0.03 &  346.67\\
                \MSR            &  0.00 &   0.03 &   0.10 &  0.00 &  20.08\\
                \LYON           &  0.03 &   1.03 &   0.88 &  0.02 &  40.55\\
                \FRANCE         &  0.09 &   3.85 &   2.47 &  0.04 &  185.13\\
            \bottomrule
        \end{tabular}
        \caption{Computing times, in seconds (OPF).}
        \label{tab:res:OPF:computing_time}
    \end{table}

    \begin{table}[!t]
        \centering
        \begin{tabular}{lrrrrr}
            \toprule
            \textbf{Solving time (secs)} &  DC & LLOA$_1$ &   LLQCP &  LLLF &  AC\\
            \midrule
                \PegaseSmall    &   0.02 &   0.43 &    1.38 &   0.03 &   55.24\\
                \PegaseMedium   &   0.09 &   1.24 &    5.35 &   0.08 &   121.96\\
                \SDET           &   0.85 &   7.85 &   79.56 &   0.85 &   124.43\\
                \RTE            &   0.06 &   3.24 &   40.07 &   0.07 &   540.92\\
                \PegaseLarge    &   0.54 &   6.70 &   52.98 &   0.53 &   675.15\\
                \PegaseHuge     &   3.70 &  15.04 &  467.52 &   2.26 &   782.49\\
                \GOC            &   0.19 &  26.98 &  272.10 &   0.24 &   1247.80\\
                \MSR            &   0.00 &   0.07 &    0.13 &   0.01 &   31.83\\
                \LYON           &   0.12 &   1.94 &    7.13 &   0.11 &   137.24\\
                \FRANCE         &   0.20 &   4.89 &   18.51 &   0.50 &   996.21\\
            \bottomrule
        \end{tabular}
        \caption{Computing times, in seconds (SCED).}
        \label{tab:res:SCED:computing_time}
    \end{table}

    First, computing times for SCED are consistently higher than for OPF.
    This is due to (i) the SCED formulation comprising more variables and constraints than OPF, and (ii) the presence of post-contingency constraints that are added lazily, thereby requiring several re-optimizations.
    The latter is especially detrimental to LLQCP: the presence of nonlinear constraints forces the use of Gurobi's barrier algorithm, which cannot be warm-started.
    Indeed, although LLQCP is surprisingly fast on the simpler OPF problems, its computing time for SCED increases by as much as 70x on the \PegaseHuge\ system.
    In contrast, linear formulations like LLLF and LLOA\textsubscript{1} make use of the simplex algorithm's warm-starting capabilities, and therefore see a lesser performance degradation between OPF and SCED.
    
    Second, for both OPF and SCED problems, LLLF is the fastest method.
    Namely, on SCED problems, LLLF outperforms LLQCP, by up to two orders of magnitude, and LLOA\textsubscript{1} by about one order of magnitude on average.
    Indeed, as mentioned in Section \ref{sec:method:discussion}, thermal constraints in LLLF are added lazily, whereas LLQCP and LLOA require the full set of transmission flow variables and constraints.
    Because, in practice, only a subset of thermal limits are binding, the LLLF formulation remains compact and requires few iterations to converge.
    Nevertheless, all SCED instances are solved by LLOA\textsubscript{1} within 30 seconds, fast enough to meet the requirements of real-time operations.
    
    Third, the performance of LLLF is not always consistent across the test cases.
    For instance, LLLF takes $0.85$s to solve the $\SDET$ system, but is 12x faster to solve the larger $\RTE$ system.
    In contrast, computing time for LLOA scales roughly linearly with the system size.
    A deeper analysis of this behavior reveals that the performance of LLLF is directly impacted by how congested a system is.
    Again, this is explained by the PTDF-based formulation used by LLLF: congested systems have more binding limits, which leads to denser formulations and slows down convergence.
    On the other hand, LLOA and LLQCP include all transmission constraints by design, and are therefore less impacted by congestion.

    

\subsection{Sensitivity analysis}
\label{sec:res:sensitivity_analysis}

\begin{figure}[!t]
        \centering
        \includegraphics[width=\linewidth]{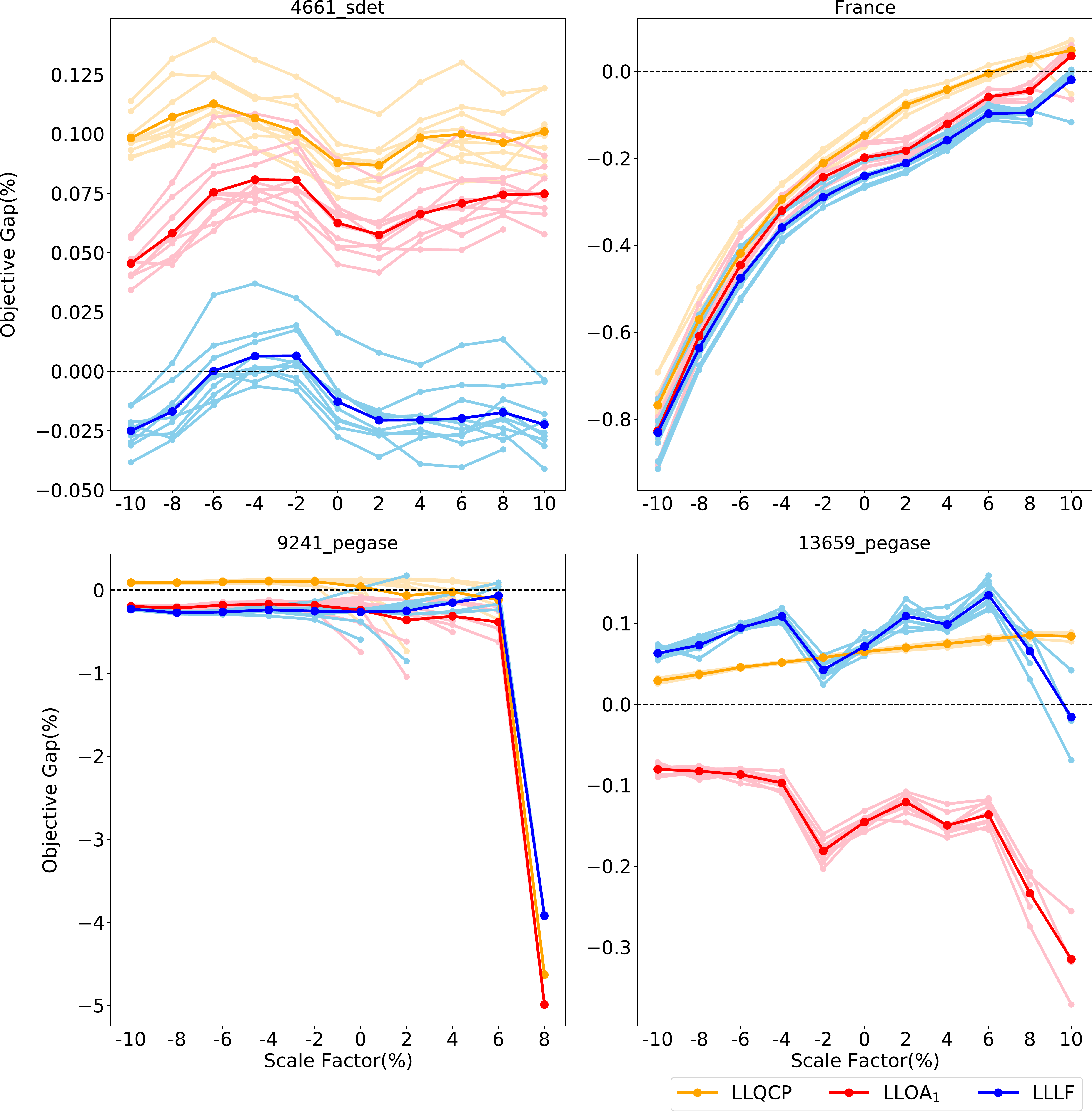}
        \caption{Behavior of the percent objective gap with respect to variations in load (SCED).}
        \label{fig:res:sensitivity:SCED:gap}
    \end{figure}

\begin{figure}[!t]
        \centering
        \includegraphics[width=\linewidth]{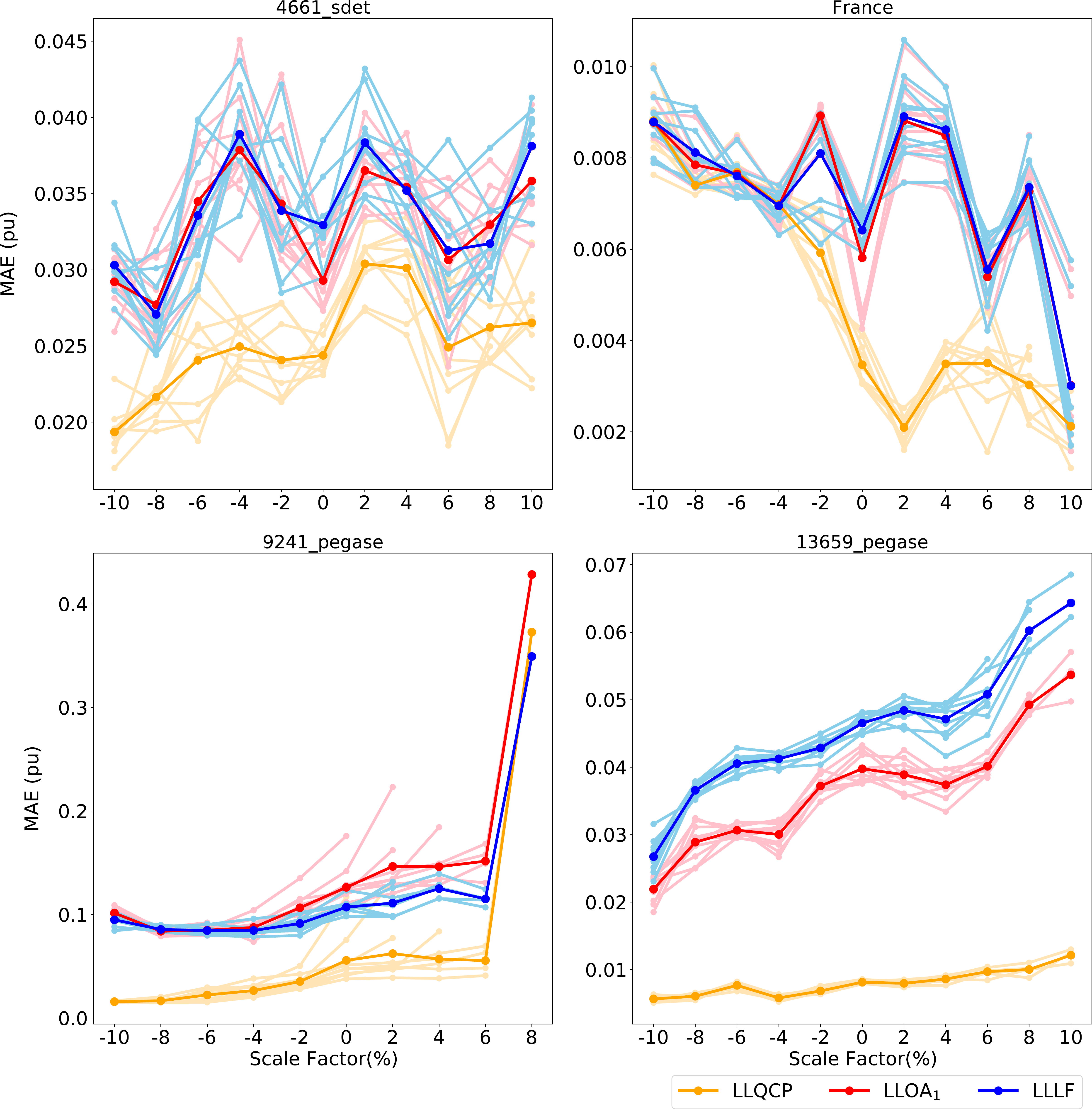}
        \caption{Behavior of the MAE with respect to variations in load (SCED).}
        \label{fig:res:sensitivity:SCED:MAE}
    \end{figure}

\begin{figure}[!t]
        \centering
        \includegraphics[width=\linewidth]{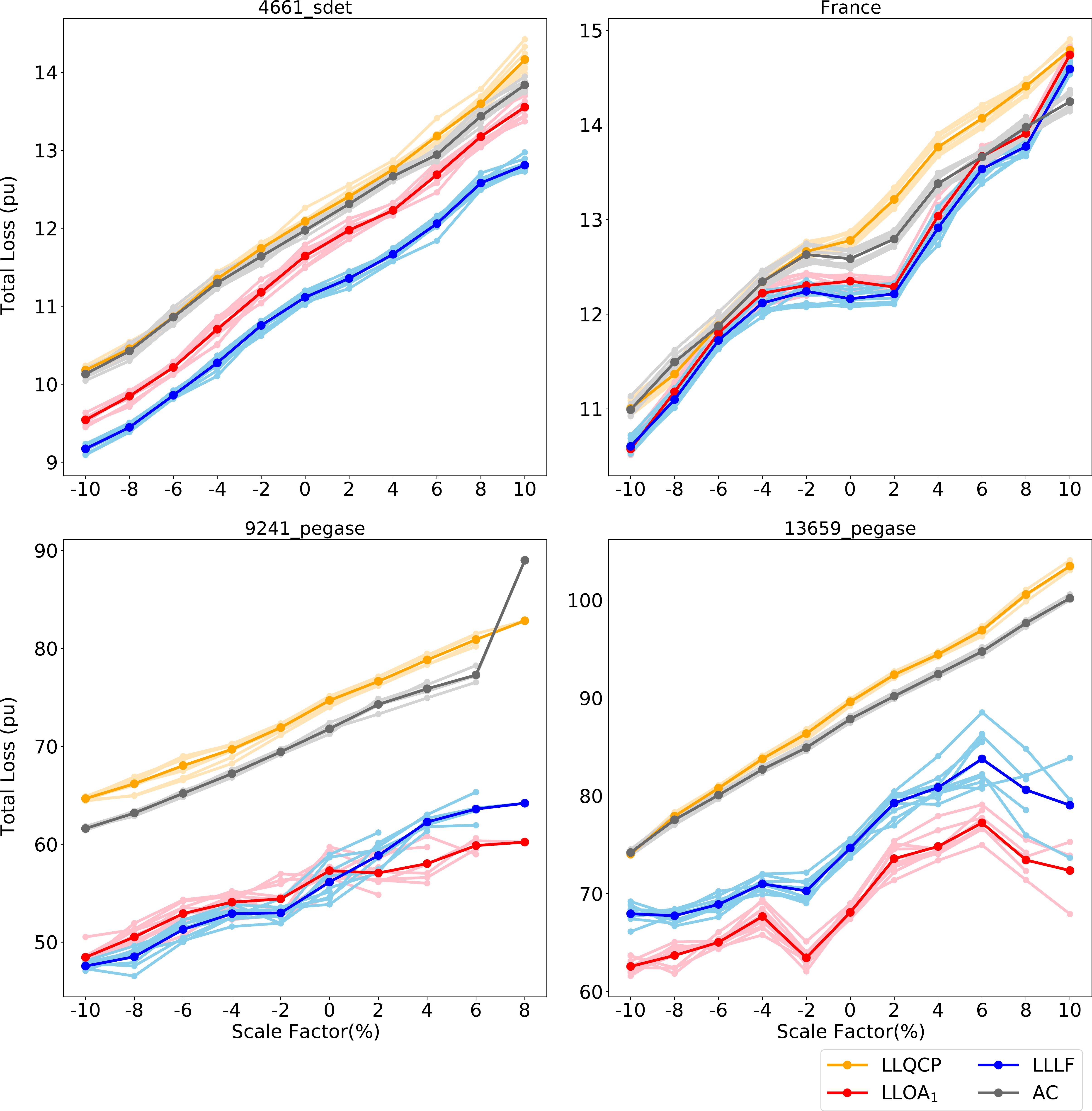}
        \caption{Behavior of the total estimated losses with respect to variations in load (SCED)}
        \label{fig:res:sensitivity:SCED:losses}
    \end{figure}

\begin{figure}[!t]
        \centering
        \includegraphics[width=\linewidth]{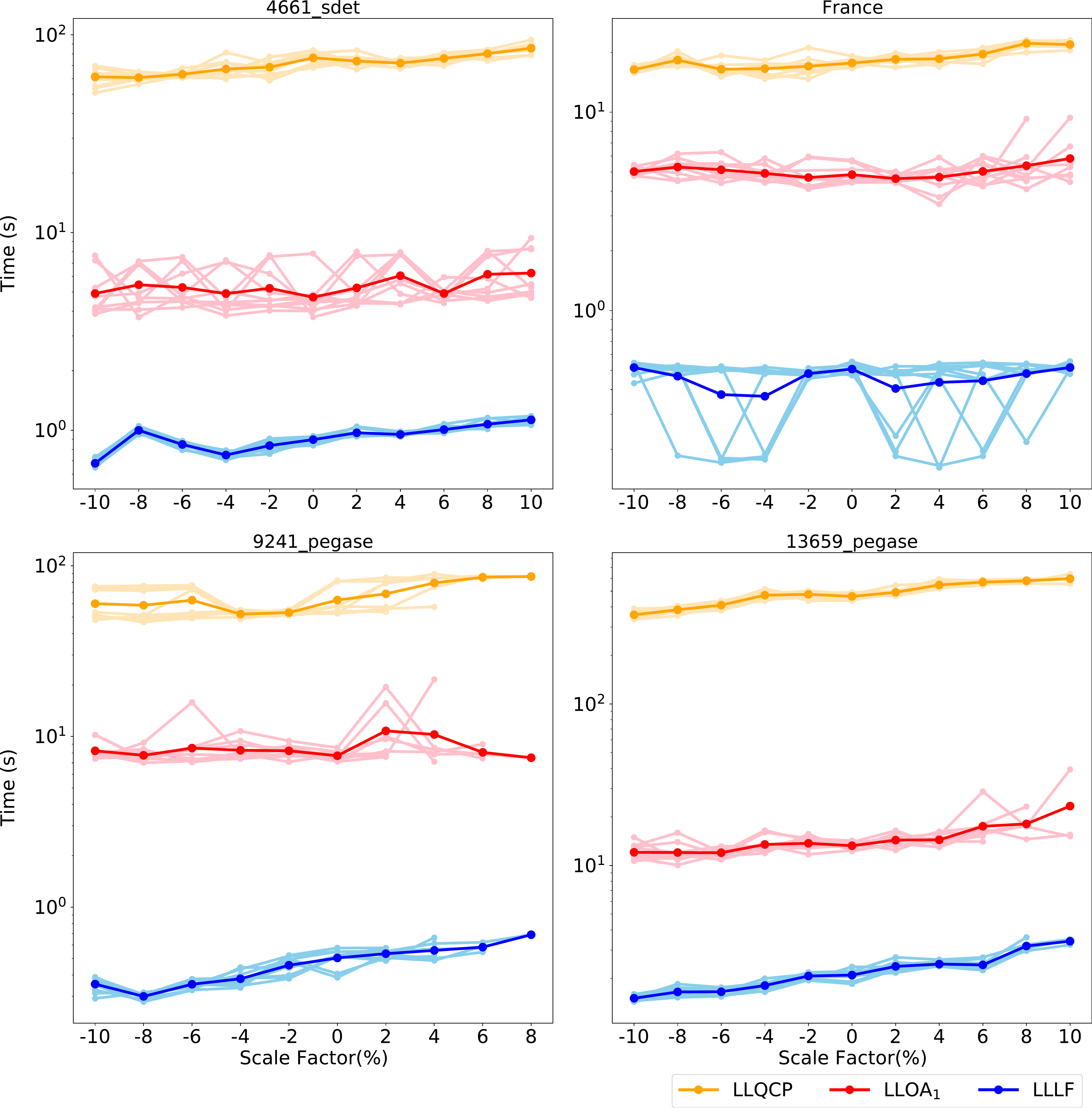}
        \caption{Behavior of the computing time with respect to variations in load (SCED).}
        \label{fig:res:sensitivity:SCED:time}
    \end{figure}

    Motivated by the above remarks, further experiments are conducted to assess the robustness of each method to changes in total load.
    Namely, additional instances are generated by perturbing each load with a multiplicative noise $\xi$, which follows a log-normal distribution with mean $\alpha$ and standard deviation $0.05$.
    In the present setting, the parameter $\alpha$, referred to as \emph{load scaling factor}, takes values between $0.9$ and $1.1$ in $0.02$ increments, corresponding to mean variations in total load that range between $-10\%$ and $+10\%$ of the reference load.
    For each value of $\alpha$, noise is sampled using ten different random seeds, which yields an extended test set of $1100$ OPF and $1100$ SCED instances.
    Numerical results for SCED are illustrated in Figures \ref{fig:res:sensitivity:SCED:gap}-\ref{fig:res:sensitivity:SCED:time}, which depict the behavior of LLLF, LLQCP, and LLOA in terms of objective gap, total losses, MAE and computing time, respectively.
    Each figure depicts results for individual seeds (in lighter tone), as well as averages across the 10 seeds.
    Instances found to be DC-infeasible, or for which Ipopt failed to find any AC-feasible solution have been factored out of the plots.
    For brevity, four representative test systems are selected: \SDET, \FRANCE, \PegaseLarge\ and \PegaseHuge; complete results for OPF and SCED are reported in the appendix.
    
    Overall, the results are consistent with those of Section \ref{sec:res:accuracy}.
    Namely, no individual method is uniformly better than the other two in terms of objective gap, LLQCP is typically more accurate in terms of MAE and total estimated losses, and LLLF displays lower computing times than LLOA, which itself outperforms LLQCP by about an order of magnitude.
    Note that the sudden changes in objective gap and MAE observed on \PegaseLarge, as total load is increased by $8\%$, are caused by a phase-transition behavior of the AC formulation, which was studied in \cite{mak2018phase}.
    
    Figures \ref{fig:res:sensitivity:SCED:gap}-\ref{fig:res:sensitivity:SCED:time} also highlight different behaviors as the load changes, depending on both the method and the test case.
    First, the behavior of LLQCP is generally smoother than LLOA and LLLF, especially on the \PegaseHuge\ system.
    This can be attributed to LLQCP being a nonlinear formulation, which is typically less sensitive to small changes in input data than linear formulations like LLOA and LLLF.
    Second, although LLLF is the fastest method, its computational efficiency is more severely affected by increases in total load than LLOA and LLQCP.
    Indeed, on SCED instances, compared to the nominal setting, LLLF is roughly 20\% faster when the load is decreased by 10\%, and 20\% slower when the load is increased by 10\%; the worst slowdown is observed on \SDET, wherein instances with the highest load are on average $2.10$x slower than in the nominal case.
    In contrast performance variations for LLOA are, on average, twice smaller than LLLF.
    As discussed in Section \ref{sec:method:discussion}, because LLLF is based on PTDF, higher congestion has an adverse impact on its performance.
    
\begin{figure}[!t]
        \centering
        \includegraphics[width=\columnwidth]{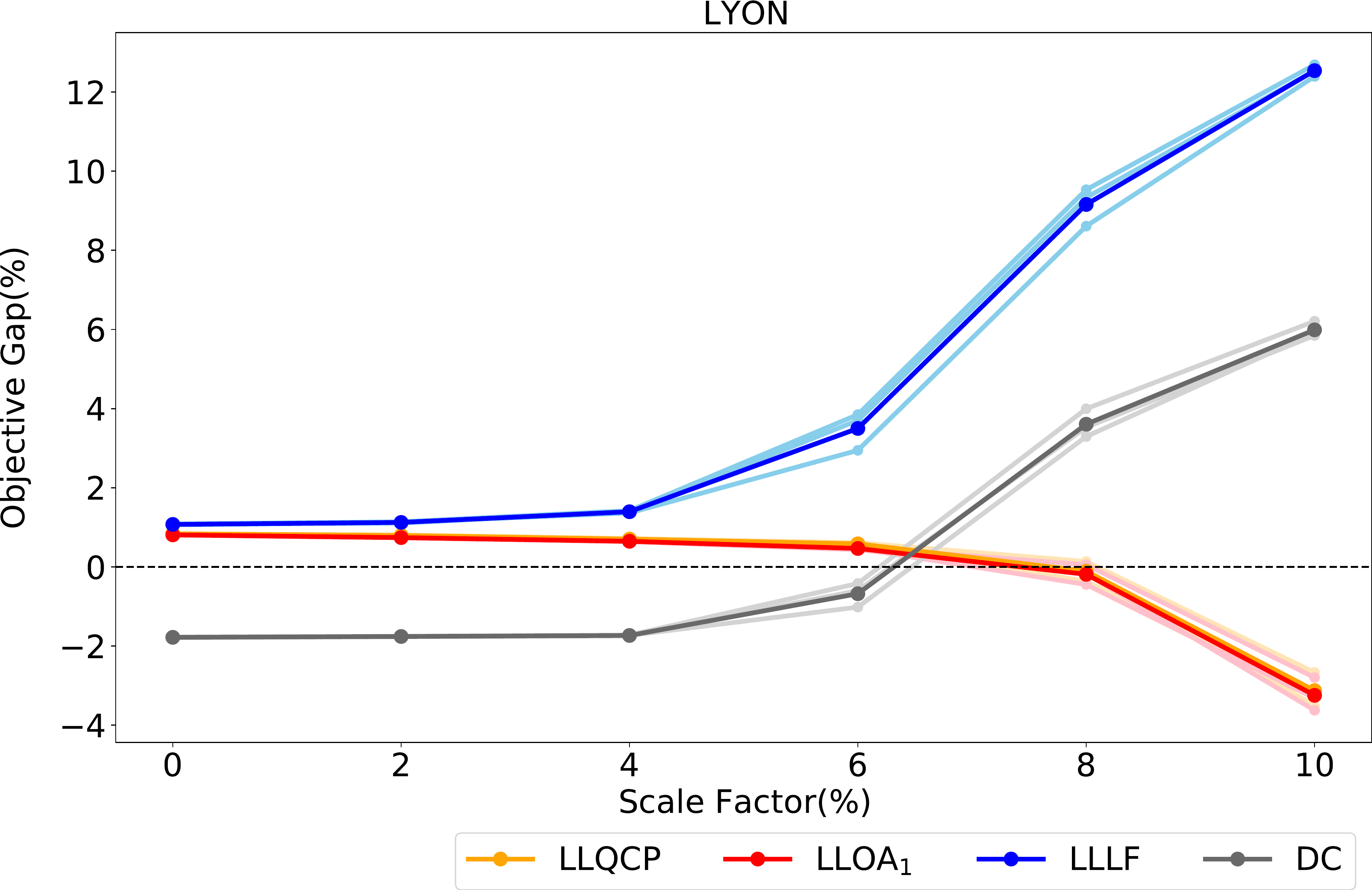}
        \caption{Objective Gap after Load Scaling for \LYON, with scale factors ranging from 0\% to 10\% and three different seeds.}
        \label{fig:load_scaling_deteriorate}
    \end{figure}

    Finally, LLOA and LLLF often behave similarly to each other in terms of objective gap, MAE and total estimated losses, most likely because they use the same reference point for estimating line losses.
    This observation holds for most of the considered test cases, OPF and SCED test cases.
    Nevertheless, a few instances show the two methods diverging as total load is increased.
    To analyze this behavior in more detail, Figure \ref{fig:load_scaling_deteriorate} reports the objective gap for each method on the \LYON\ SCED; results from only 3 seeds are depicted for readability.
    While LLOA and LLQCP are almost indistinguishable, with objective gap that remains low, the objective gap of LLLF increases to over 12\%, meaning the objective value of LLLF is significantly higher than LLOA and LLQCP.
    The observed behavior remains even when using the AC solution as a reference point for computing the loss factors. This behaviour is not an exception on the \LYON\ system, we observed similar pattern on $\texttt{500\_goc}$ as well.
    A deeper analysis reveals that the change in objective value is caused by several thermal violations in LLLF, which trigger higher reserve levels and transmission violation penalties; the same set transmission lines is at fault across the different seeds. 
    Furthermore, although the thermal limits are not violated in either LLQCP, LLOA, or AC formulations, taking the corresponding solutions and evaluating the flows as per Eq. \eqref{eq:LLLF:thermal_limits} systematically yields violations on those branches.
    This demonstrates the limitation of LLLF in estimating line flows accurately, and that such inaccuracies can result in significantly higher operating costs.

    
%% file: tex/power_flow.tex
\section{AC Feasibility Restoration}
\label{sec:power_flow}

While electricity markets are based on DC formulations, operators must ensure that final dispatch instructions are AC-feasible.
To that end, the market solution seeds a power flow analysis to determine voltage levels and active/reactive power dispatch for each generator.
Whenever the power flow engine is unable to recover a feasible AC operating point, operators must manually override some of the market decisions, which typically leads to increased production costs.

This process is replicated by seeding a power flow engine with the solution of DC, LLQCP, LLOA, and LLLF methods.
Although the power flow engine does not enforce operating and thermal limits, violations in reactive power capacity are alleviated --but not always eliminated-- via PV-PQ bus switching.
Recall that reserve dispatch (in SCED) and total operating costs only depend on active power dispatch which, by definition, the power flow engine does not modify except for the generator located at the slack bus.
Therefore, the objective value will not be changed much by the power flow and, accordingly, the rest of the analysis will focus on AC feasibility.
To that end, Tables \ref{tab:res:PowerFlow:OPF:PQ-Viol}-\ref{tab:res:PowerFlow:OPF:VT-Viol} report, for the OPF instances, the number (\#viol) and largest magnitude (Max) of constraint violations for the power-power flow solutions, namely, active/reactive power limits, voltage limits, and thermal limits; results for SCED instances are in Tables \ref{tab:res:PowerFlow:SCED:PQ-Viol}-\ref{tab:res:PowerFlow:SCED:VT-Viol}.
All results are obtained in the nominal setting, i.e., without any load perturbation, and test cases for which the power flow did not converge are factored out of the tables.

Results are consistent across OPF and SCED instances.
Vanilla DC leads to the highest violations after power flow, especially regarding active power capacity and thermal limits.
The former is due to the slack generator having to compensate for all losses, which are ignored in vanilla DC yet often exceed that generator's maximum output.
Similarly, large thermal limit violations are found to systematically occur on lines adjacent to the slack bus.
In contrast, LLLF, LLQCP, and LLOA exhibit similar behavior, with fewer and smaller violations than vanilla DC, although some violations still occur.
In particular, active power capacity violations are directly correlated to the accuracy of total estimated losses: the more accurate the estimated losses, the smaller the violation.
The use of a more robust power flow engine --which is beyond the scope of this paper-- would likely reduce overall violations in reactive power capacity.
Nevertheless, the observed violations in reactive power and voltage magnitude highlight the limitations of DC-based formulations, which ignore reactive power and assume voltage magnitude to be 1 per-unit.

%% file: tex/appendix.tex
\appendix

%
%
\begin{table}[H]
    \centering
    \begin{tabular}{llrrrr}
        \toprule
            && \multicolumn{2}{c}{Active Power} & \multicolumn{2}{c}{Reactive Power}\\
            \cmidrule(lr){3-4} \cmidrule(lr){5-6}
            Test case & Method & \#viol. & Max & \#viol. & Max\\
        \midrule
        \multirow{4}{*}{\PegaseSmall} & DC & 1 & 19.21  & 0 & 0.00    \\
                & LLQCP & 1 & 1.62   & 0 & 0.00    \\
                & LLOA$_{1}$ & 1 & 1.80   & 0 & 0.00    \\
                & LLLF & 1 & 2.42   & 0 & 0.00    \\
        \midrule
        \multirow{4}{*}{\PegaseMedium} & DC & 1 & 39.69  & 0 & 0.00 \\
                & LLQCP & 1 & 3.06   & 0 & 0.00    \\
                & LLOA$_{1}$ & 1 & 13.56  & 0 & 0.00    \\
                & LLLF & 1 & 14.46  & 0 & 0.00    \\
        \midrule
        \multirow{4}{*}{\SDET} & DC & 0 & 0.00   & 0 & 0.00 \\
                & LLQCP & 0 & 0.00   & 0 & 0.00 \\
                & LLOA$_{1}$ & 0 & 0.00   & 0 & 0.00 \\
                & LLLF & 0 & 0.00   & 0 & 0.00 \\
        \midrule
        \multirow{4}{*}{\MSR} & DC & 0 & 0.00   & 2 & 1.86 \\
                & LLQCP & 0 & 0.00      & 1 & 1.87 \\
                & LLOA$_{1}$ & 0 & 0.00   & 1 & 1.87 \\
                & LLLF & 0 & 0.00   & 1 & 1.87 \\
        \midrule
        \multirow{4}{*}{\LYON} & DC & 1 & 12.14  & 6 & 1.13 \\
                & LLQCP & 0 & 0.00   & 6 & 0.83 \\
                & LLOA$_{1}$ & 0 & 0.00   & 6 & 0.83 \\
                & LLLF & 1 & 0.65   & 6 & 0.83 \\
        \midrule
        \multirow{4}{*}{\FRANCE} & DC & 1 & 14.95  & 5 & 1.56 \\
                & LLQCP & 1 & 0.72   & 4 & 1.56 \\
                & LLOA$_{1}$ & 1 & 1.00   & 5 & 1.56 \\
                & LLLF & 1 & 1.19   & 5 & 1.55 \\
        \bottomrule
    \end{tabular}
    \caption{Constraint violation statistics for active and reactive capacity (OPF). Violations are in per-unit.}
    \label{tab:res:PowerFlow:OPF:PQ-Viol}
\end{table}

\begin{table}[H]
    \centering
    \begin{tabular}{llrrrr}
        \toprule
            && \multicolumn{2}{c}{Voltage magnitude} & \multicolumn{2}{c}{Thermal}\\
            \cmidrule(lr){3-4} \cmidrule(lr){5-6}
            Test case & Method & \#viol. & Max & \#viol. & Max\\
        \midrule
        \multirow{4}{*}{\PegaseSmall} & DC & 1    & 0.01 & 23  & 6.05  \\
                & LLQCP & 0    & 0.00    & 21  & 4.06  \\
                & LLOA$_{1}$ & 0    & 0.00    & 21  & 4.05  \\
                & LLLF & 0    & 0.00    & 22  & 4.35  \\
        \midrule
        \multirow{4}{*}{\PegaseMedium} 
                & LLQCP & 90   & 0.04 & 33  & 5.70  \\
                & LLOA$_{1}$ & 49   & 0.03 & 32  & 5.64  \\
                & LLLF & 115  & 0.05 & 34  & 6.52  \\
        \midrule
        \multirow{4}{*}{\SDET} & DC & 9    & 0.02 & 196 & 5.30  \\
                & LLQCP & 8    & 0.02 & 180 & 1.83  \\
                & LLOA$_{1}$ & 8    & 0.02 & 177 & 2.13  \\
                & LLLF & 8    & 0.02 & 190 & 2.24  \\
        \midrule
        \multirow{4}{*}{\MSR} & DC & 36   & 0.02 & 0   & 0.00     \\
                & LLQCP & 35   & 0.02 & 0   & 0.00     \\
                & LLOA$_{1}$ & 35   & 0.02 & 0   & 0.00     \\
                & LLLF & 35   & 0.02 & 0   & 0.00     \\
        \midrule
        \multirow{4}{*}{\LYON} & DC & 542  & 0.06 & 7   & 13.07 \\
                & LLQCP & 538  & 0.06 & 1   & 0.11  \\
                & LLOA$_{1}$ & 539  & 0.06 & 1   & 0.11  \\
                & LLLF & 538  & 0.06 & 3   & 0.12  \\
        \midrule
        \multirow{4}{*}{\FRANCE} & DC & 413  & 0.05 & 10  & 18.37 \\
                & LLQCP & 409  & 0.05 & 9   & 0.18  \\
                & LLOA$_{1}$ & 403  & 0.05 & 9   & 0.18  \\
                & LLLF & 404  & 0.05 & 9   & 0.21  \\
        \bottomrule
    \end{tabular}
    \caption{Constraint violation statistics for voltage magnitude and thermal limits (OPF). Violations are in per-unit.}
    \label{tab:res:PowerFlow:OPF:VT-Viol}
\end{table}

\begin{table}[!t]
    \centering
    \begin{tabular}{llrrrr}
        \toprule
            && \multicolumn{2}{c}{Active Power} & \multicolumn{2}{c}{Reactive Power}\\
            \cmidrule(lr){3-4} \cmidrule(lr){5-6}
            Test case & Method & \#viol. & Max & \#viol. & Max\\
        \midrule
        \multirow{4}{*}{\PegaseSmall} & DC & 1 & 18.15 & 0  & 0.00    \\
                & LLQCP & 1 & 1.63  & 0  & 0.00    \\
                & LLOA$_{1}$ & 1 & 1.68  & 0  & 0.00    \\
                & LLLF & 1 & 2.29  & 0  & 0.00    \\
        \midrule
        \multirow{4}{*}{\PegaseMedium} & DC & 1 & 35.79 & 0  & 0.00    \\
                & LLQCP & 1 & 2.75  & 0  & 0.00    \\
                & LLOA$_{1}$ & 1 & 5.13  & 0  & 0.00    \\
                & LLLF & 1 & 6.89  & 0  & 0.00    \\
        \midrule
        \multirow{4}{*}{\SDET} & DC & 0 & 0.00  & 0  & 0.00 \\
                & LLQCP & 0 & 0.00  & 0  & 0.00 \\
                & LLOA$_{1}$ & 0 & 0.00  & 0  & 0.00 \\
                & LLLF & 0 & 0.00  & 0  & 0.00 \\
        \midrule
        \multirow{4}{*}{\MSR} & DC & 0 & 0.00  & 7  & 1.43 \\
                & LLQCP & 0 & 0.00     & 7  & 1.45 \\
                & LLOA$_{1}$ & 0 & 0.00  & 7  & 1.45 \\
                & LLLF & 0 & 0.00  & 7  & 1.45 \\
        \midrule
        \multirow{4}{*}{\LYON} & DC & 1 & 11.94 & 19 & 1.42 \\
                & LLQCP & 0 & 0.00  & 17 & 1.42 \\
                & LLOA$_{1}$ & 0 & 0.00  & 17 & 1.42 \\
                & LLLF & 1 & 0.75  & 19 & 1.44 \\
        \midrule
        \multirow{4}{*}{\FRANCE} & DC & 1 & 14.28 & 19 & 1.28 \\
                & LLQCP & 1 & 0.35  & 20 & 1.28 \\
                & LLOA$_{1}$ & 1 & 0.60  & 21 & 1.28 \\
                & LLLF & 1 & 0.77  & 21 & 1.28 \\
        \bottomrule
    \end{tabular}
    \caption{Active and reactive capacity violation statistics (SCED). Constraint violations are in per-unit.}
    \label{tab:res:PowerFlow:SCED:PQ-Viol}
\end{table}

\begin{table}[!t]
    \centering
    \begin{tabular}{llrrrr}
        \toprule
            && \multicolumn{2}{c}{Voltage magnitude} & \multicolumn{2}{c}{Thermal}\\
            \cmidrule(lr){3-4} \cmidrule(lr){5-6}
            Test case & Method & \#viol. & Max & \#viol. & Max\\
        \midrule
        \multirow{4}{*}{\PegaseSmall} & DC & 1    & 0.02 & 22  & 5.32  \\
                & LLQCP & 1    & 0.01 & 19  & 3.41  \\
                & LLOA$_{1}$ & 1    & 0.01 & 19  & 3.41  \\
                & LLLF & 1    & 0.01 & 18  & 3.39  \\
        \midrule
        \multirow{4}{*}{\PegaseMedium} & DC & 73   & 0.09 & 53  & 9.72  \\
                & LLQCP & 18   & 0.02 & 35  & 5.42  \\
                & LLOA$_{1}$ & 18   & 0.02 & 33  & 5.34  \\
                & LLLF & 16   & 0.02 & 34  & 5.71  \\
        \midrule
        \multirow{4}{*}{\SDET} & DC & 9    & 0.02 & 176 & 5.15  \\
                & LLQCP & 8    & 0.02 & 166 & 1.22   \\
                & LLOA$_{1}$ & 8    & 0.02 & 170 & 1.22   \\
                & LLLF & 8    & 0.02 & 177 & 1.42   \\
        \midrule
        \multirow{4}{*}{\MSR} & DC & 70   & 0.02 & 0   & 0.00      \\
                & LLQCP & 70   & 0.02 & 0   & 0.00      \\
                & LLOA$_{1}$ & 70   & 0.02 & 0   & 0.00      \\
                & LLLF & 70   & 0.02 & 0   & 0.00      \\
        \midrule
        \multirow{4}{*}{\LYON} & DC & 573  & 0.06 & 5   & 12.63 \\
                & LLQCP & 572  & 0.06 & 4   & 0.20   \\
                & LLOA$_{1}$ & 572  & 0.06 & 4   & 0.20   \\
                & LLLF & 566  & 0.06 & 3   & 0.14   \\
        \midrule
        \multirow{4}{*}{\FRANCE} & DC & 479  & 0.05 & 6   & 17.26 \\
                & LLQCP & 467  & 0.05 & 6   & 0.17   \\
                & LLOA$_{1}$ & 467  & 0.05 & 6   & 0.17   \\
                & LLLF & 467  & 0.05 & 7   & 0.20   \\
        \bottomrule
    \end{tabular}
    \caption{Voltage magnitude and thermal violation statistics (SCED). Constraint violations are in per-unit.}
    \label{tab:res:PowerFlow:SCED:VT-Viol}
\end{table}

\begin{figure}
    \centering
    \includegraphics[width=\linewidth]{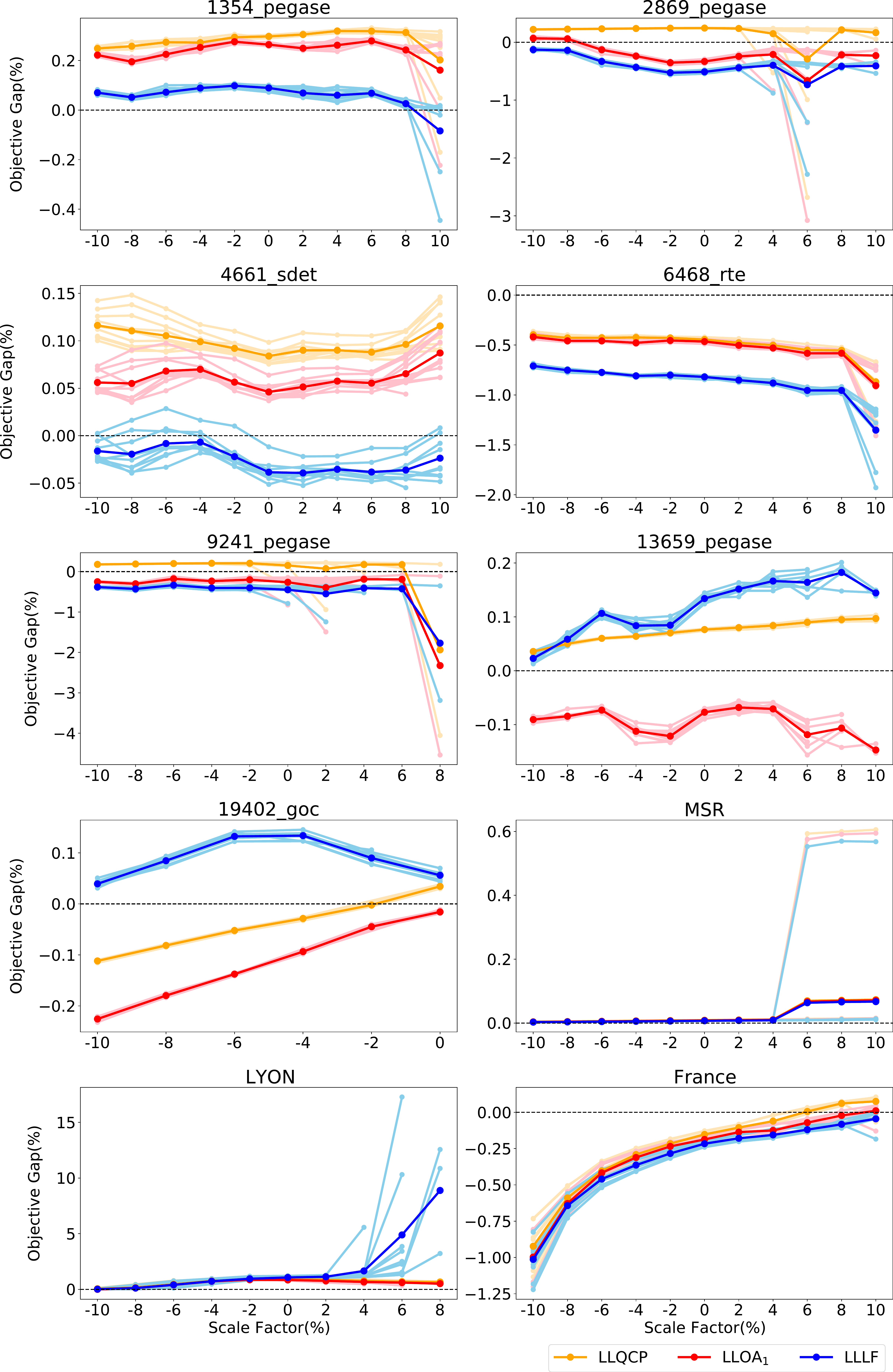}
    \caption{Behavior of the percent objective gap with respect to variations in load (OPF).}
    \label{fig:opf_load_scaling_perturbed}
\end{figure}

\begin{figure}
    \centering
    \includegraphics[width=\linewidth]{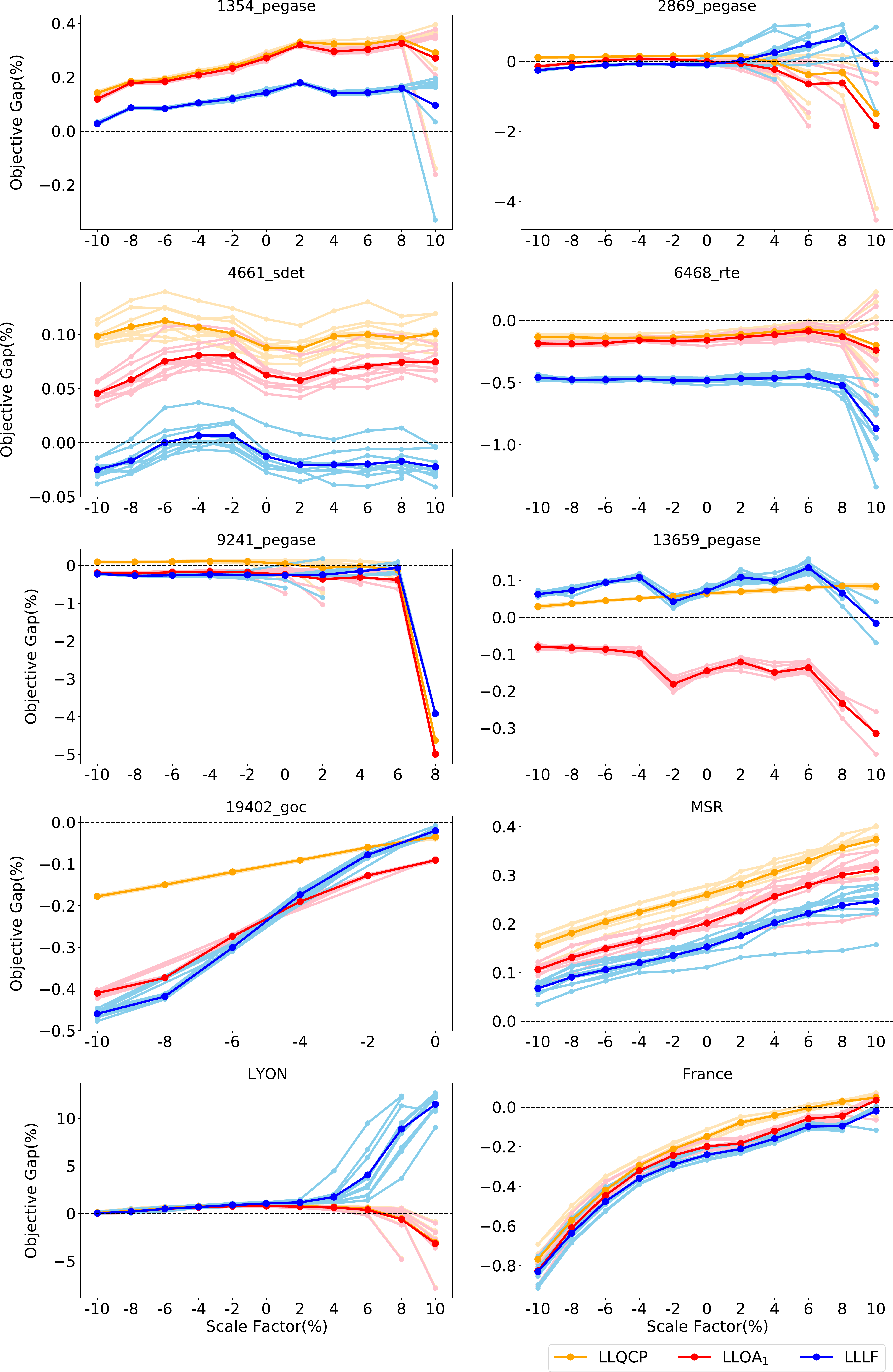}
    \caption{Behavior of the percent objective gap with respect to variations in load (SCED).}
    \label{fig:sced_load_scaling_perturbed}
\end{figure}

\begin{figure}
    \centering
    \includegraphics[width=\linewidth]{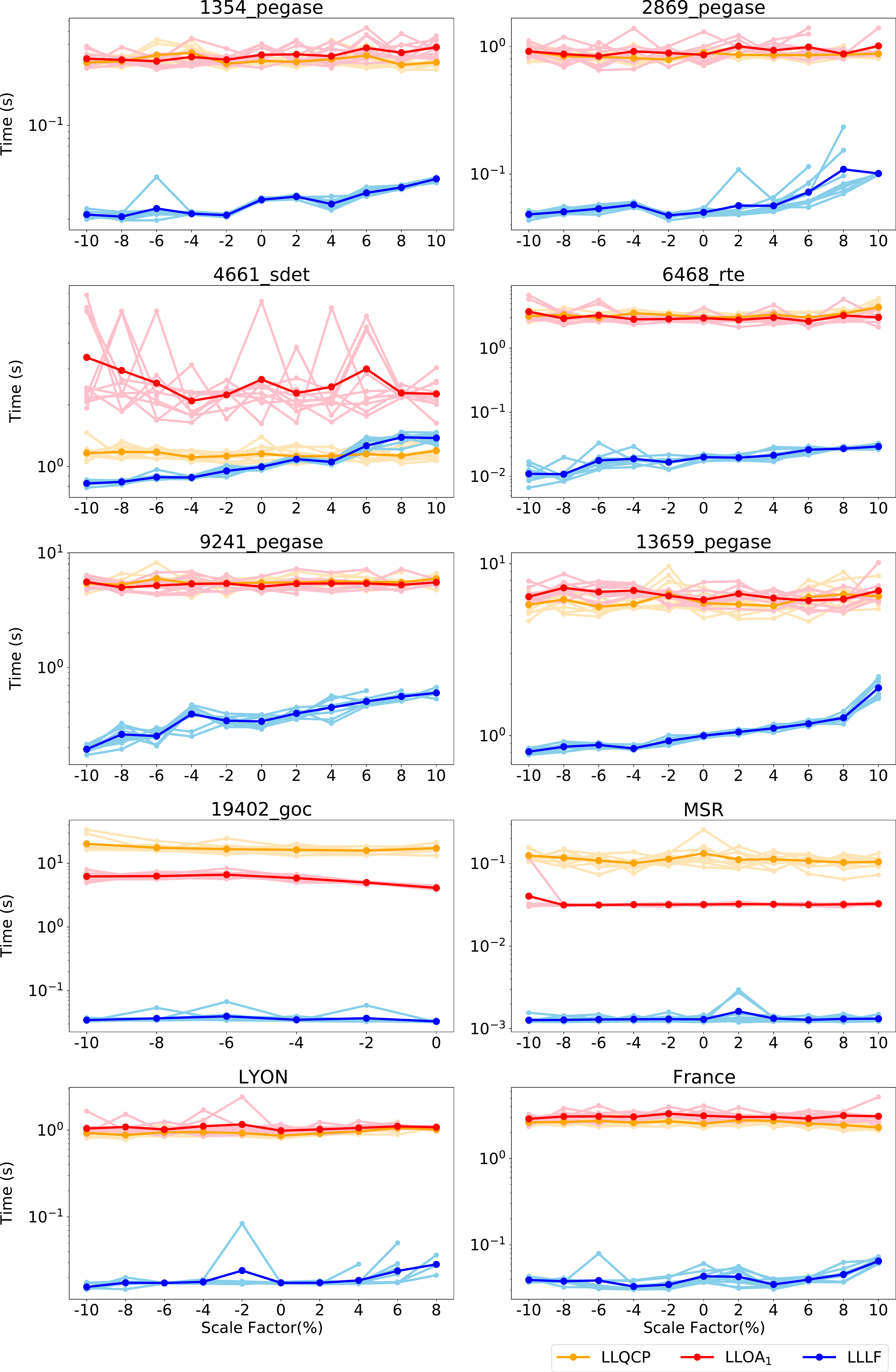}
    \caption{Behavior of the computing time with respect to variations in load (OPF).}
    \label{fig:Load_Scaling_10seed_OPF_time}
\end{figure}

\begin{figure}
    \centering
    \includegraphics[width=\linewidth]{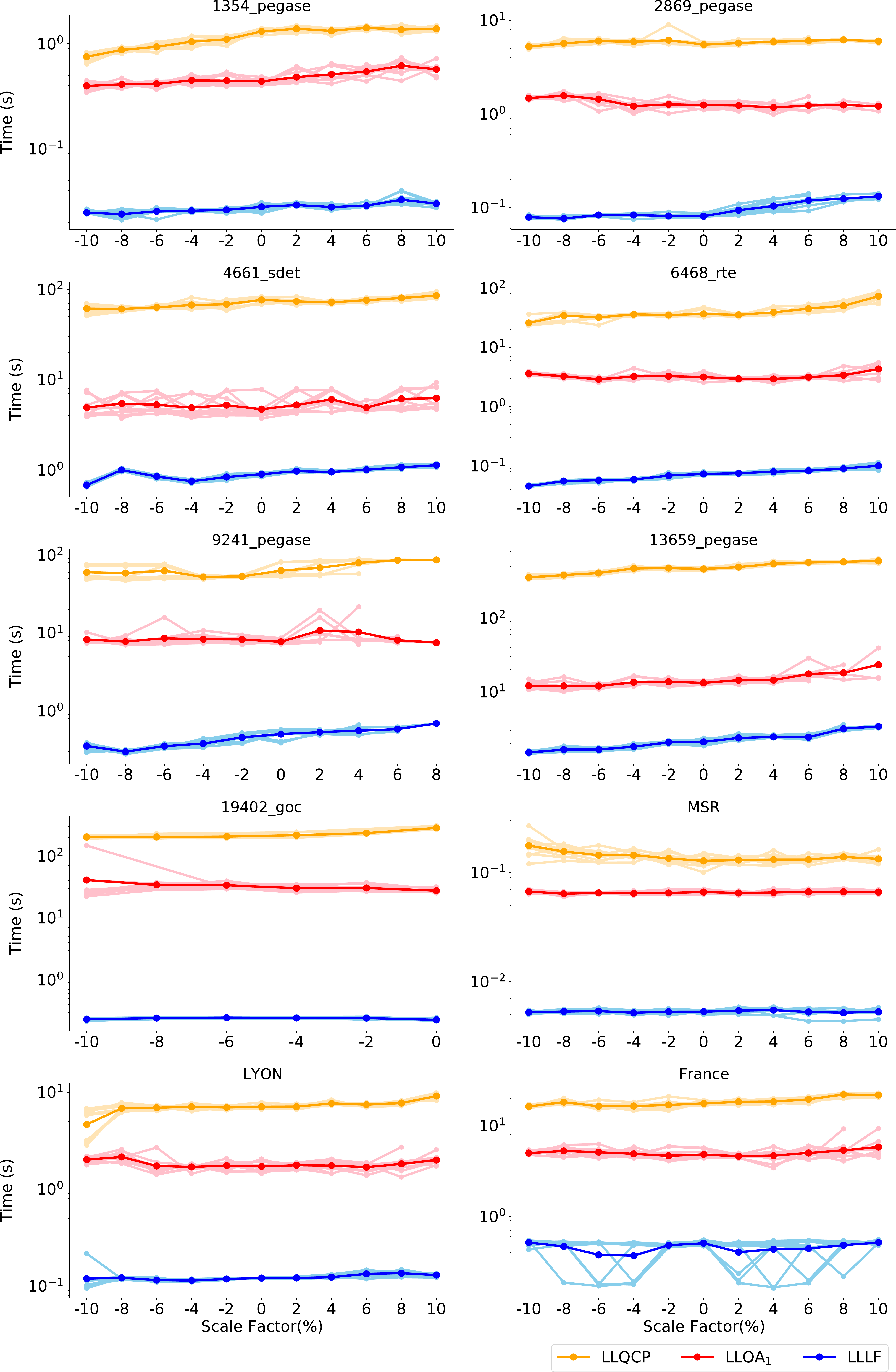}
    \caption{Behavior of the computing time with respect to variations in load (SCED).}
    \label{fig:Load_Scaling_10seed_SCED_time}
\end{figure}

\begin{figure}[!t]
    \centering
    \includegraphics[width=\columnwidth]{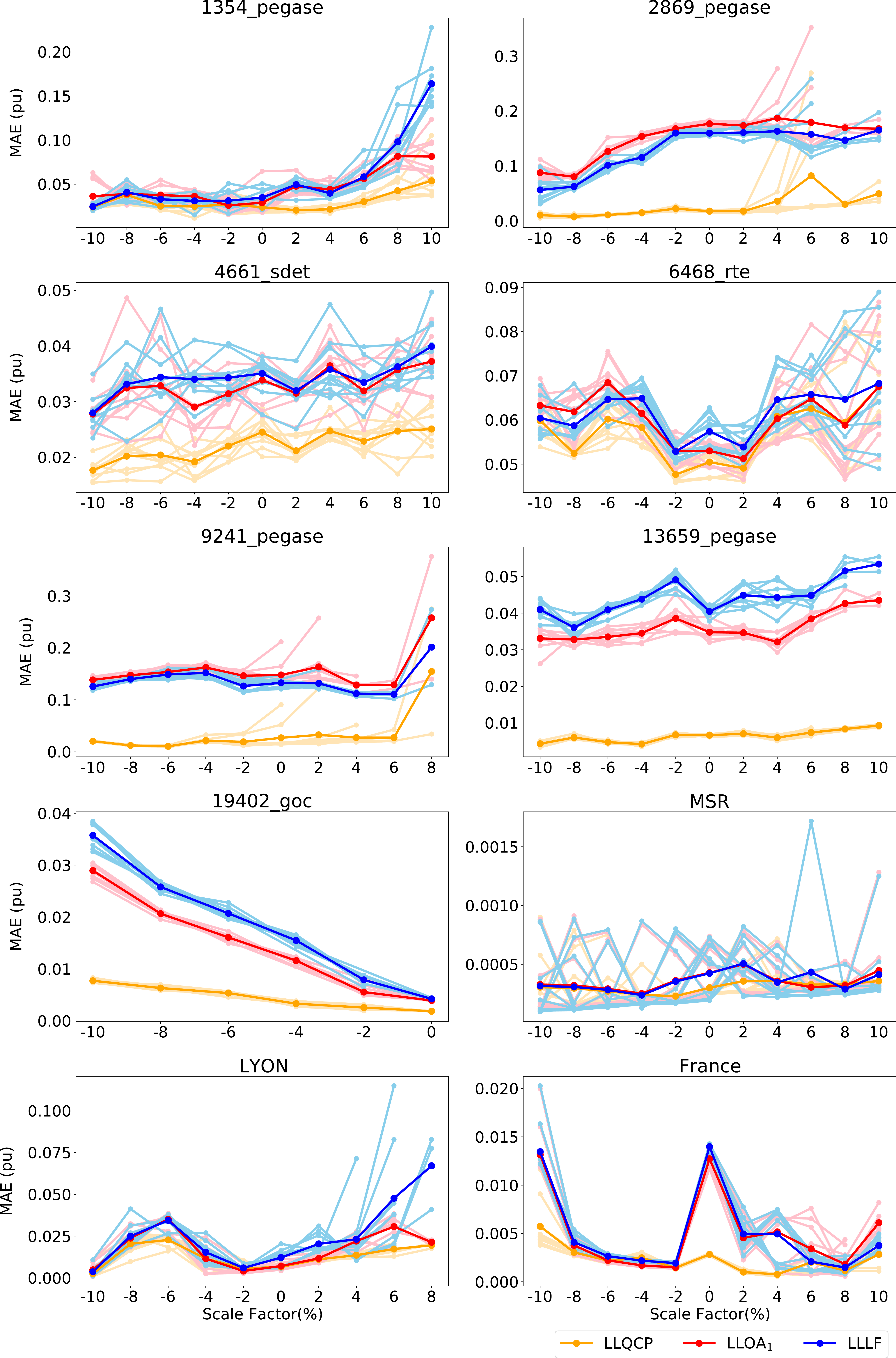}
    \caption{Behavior of the MAE with respect to variations in load (OPF).}
    \label{fig:Load_Scaling_10seed_OPF_mae}
\end{figure}

\begin{figure}[!t]
    \centering
    \includegraphics[width=\columnwidth]{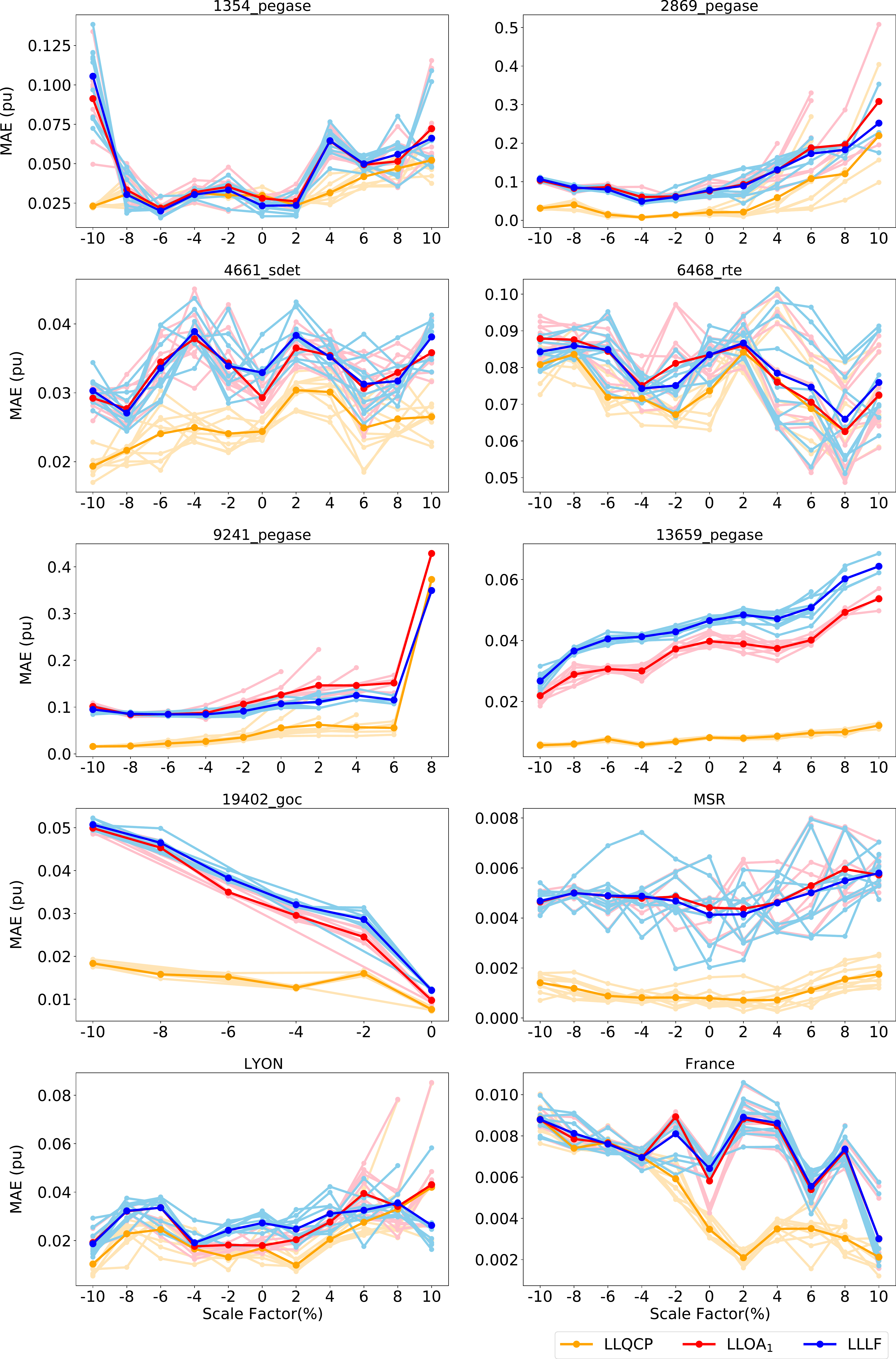}
    \caption{Behavior of the MAE with respect to variations in load (SCED).}
    \label{fig:Load_Scaling_10seed_SCED_mae}
\end{figure}

\begin{figure}
    \centering
    \includegraphics[width=\linewidth]{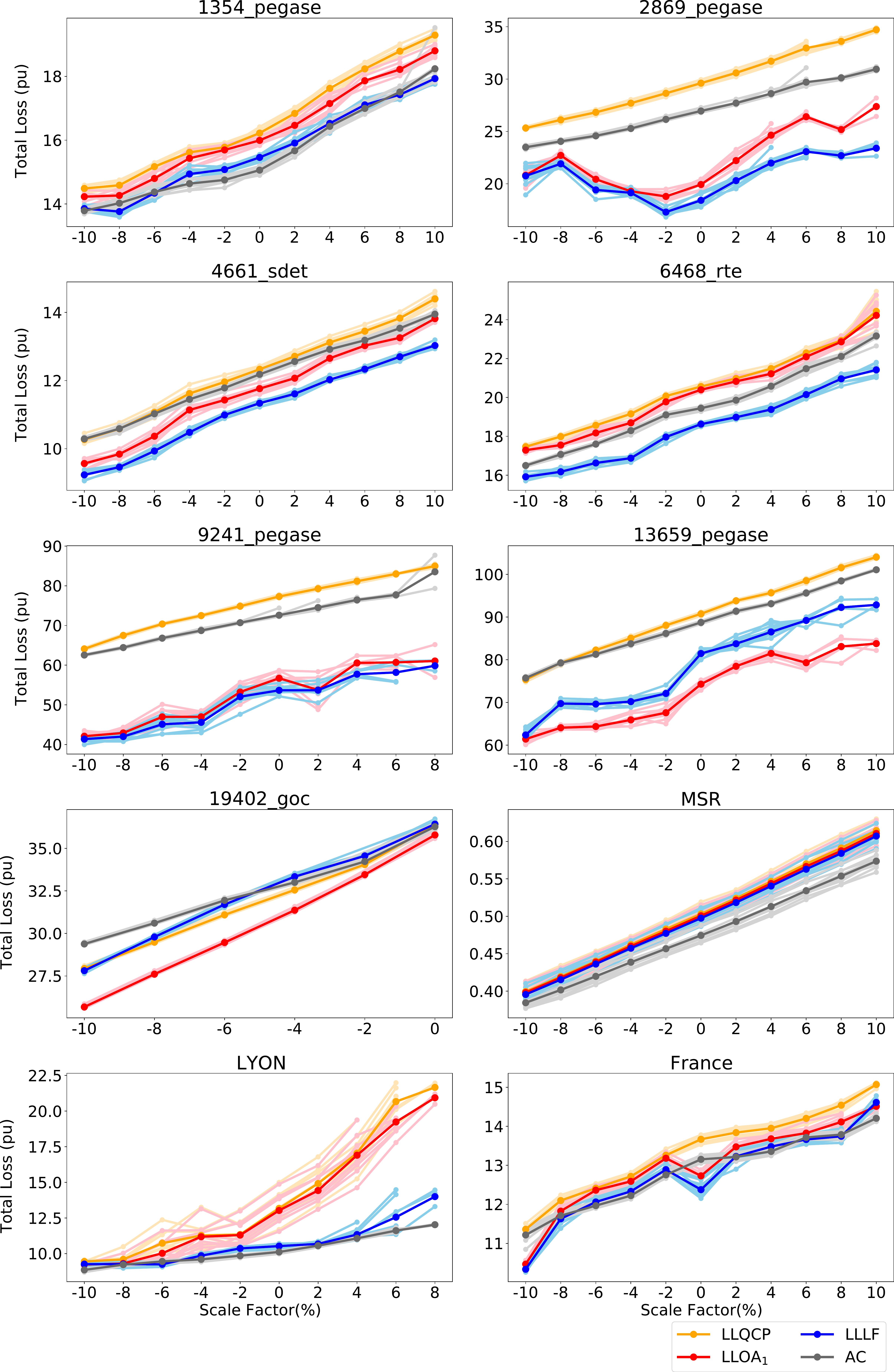}
    \caption{Behavior of the total estimated losses with respect to variations in load (OPF)}
    \label{fig:Load_Scaling_10seed_OPF_total_loss}
\end{figure}

\begin{figure}
    \centering
    \includegraphics[width=\linewidth]{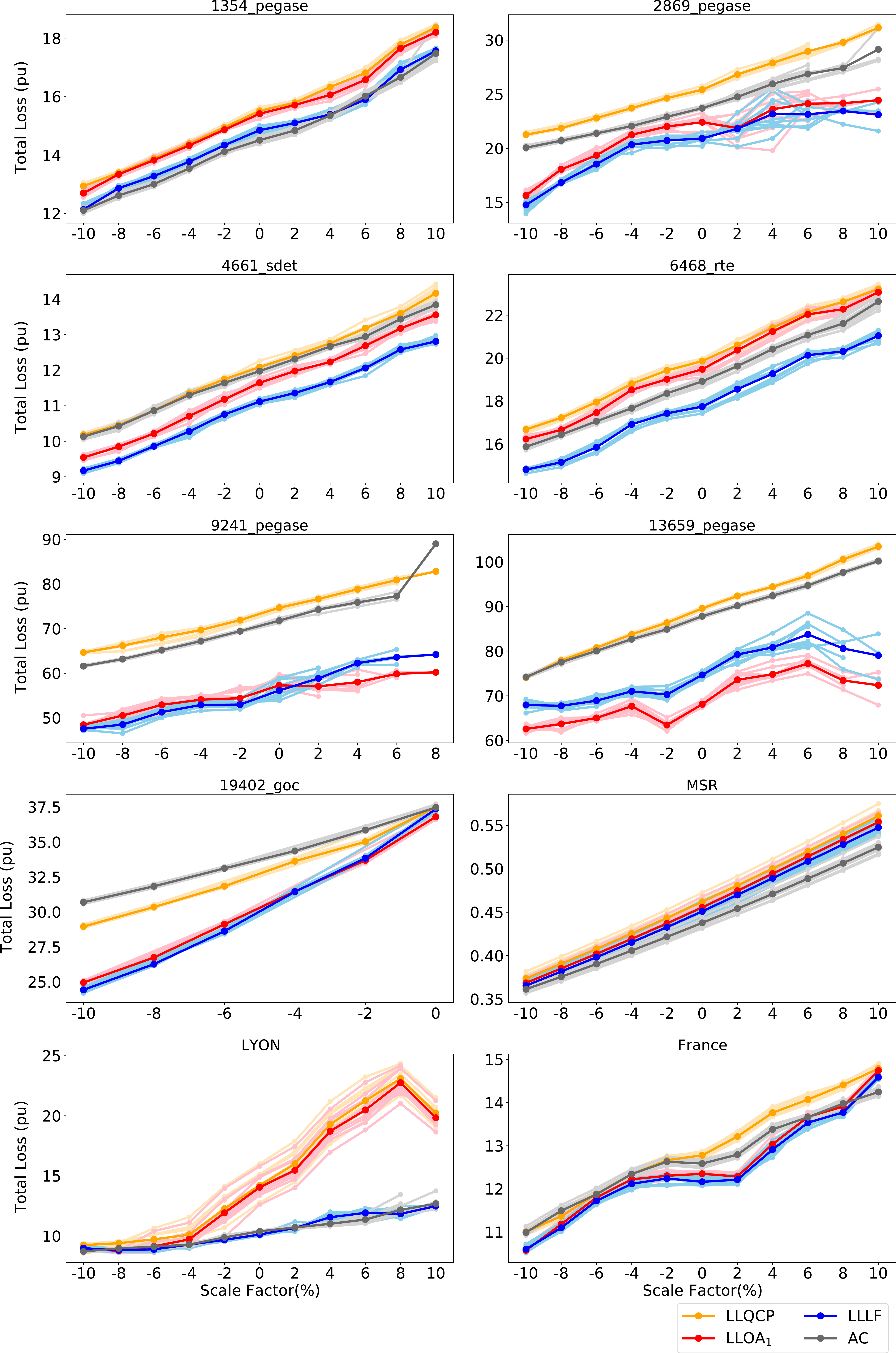}
    \caption{Behavior of the total estimated losses with respect to variations in load (SCED)}
    \label{fig:Load_Scaling_10seed_SCED_total_loss}
\end{figure}

\begin{table}[!t]
        \centering
        \begin{tabular}{lrrrr}
            \toprule
            \textbf{Obj Gap(\%)} &   DC &   LLOA$_1$ &  LLQCP &   LLLF \\
            \midrule
            \PegaseSmall &  1.66 &   \textbf{0.26} &   0.30 &   0.27 \\
            \PegaseMedium &  2.09 &   1.34 &   \textbf{0.40} &  1.04 \\
            \SDET &   1.28 &   0.49 &   \textbf{0.42} & 0.53\\
            \RTE &  6.76 &  4.16 & \textbf{4.10} &  4.41 \\
            \PegaseLarge &  2.17 &  1.24 &   \textbf{0.22} &  1.21 \\
            \PegaseHuge &  1.83 &  0.35 &  \textbf{ 0.12} &   0.29 \\
            \GOC &  -1.60 &  \textbf{0.01} &  0.03 &  0.06 \\
            \MSR &  2.45 &   0.06 &   \textbf{0.06} &   0.06 \\
            \LYON &   4.88 &   \textbf{0.83} &   0.87 &   1.06 \\
            \FRANCE &  6.71 &   0.86 &  \textbf{0.69} &  0.90 \\
            \bottomrule
        \end{tabular}
        \caption{Objective Differences in Percentage for OPF after running load flow.}
    \label{tab:OPF_load_flow_Obj_gap}
\end{table}

\begin{table}[!t]
        \centering
        \begin{tabular}{lrrrr}
            \toprule
            \textbf{Obj Gap(\%)} &   DC &   LLOA$_1$ &  LLQCP &   LLLF \\
            \midrule
            \PegaseSmall &  1.53 &   0.29 &   0.30 &   \textbf{0.16} \\
            \PegaseMedium &  1.50 &   0.28 &   \textbf{0.21} &  0.36 \\
            \SDET &   1.20 &   0.42 &   \textbf{0.37} & 0.47\\
            \RTE &  6.45 &  \textbf{4.29} & 4.30 &  4.46 \\
            \PegaseLarge &  2.41 &  1.22 &   \textbf{0.62} &  0.93 \\
            \MSR &  -3.74 &   0.26 &   0.30 &   \textbf{0.22} \\
            \LYON &   4.52 &   \textbf{0.77} &   0.79 &   1.14 \\
            \FRANCE &  6.21 &   0.92 &  \textbf{0.73} &  1.00 \\
            \bottomrule
        \end{tabular}
        \caption{Objective Differences in Percentage for SCED after running load flow.}
    \label{tab:SCED_load_flow_Obj_gap}
\end{table}

\begin{table}[!t]
        \centering
        \begin{tabular}{lrrrr}
            \toprule
            \textbf{MAE(pu)} &   DC &   LLOA$_1$ &  LLQCP &   LLLF \\
            \midrule
            \PegaseSmall &  0.06 &   0.00 &   0.00 &   0.00 \\
            \PegaseMedium &  0.06 &   0.02 &   0.01 &  0.02 \\
            \SDET &   0.02 &   0.00 &   0.00 & 0.00\\
            \RTE &  0.12 &  0.08 & 0.09 &  0.09 \\
            \PegaseLarge &  0.06 &  0.02 &   0.00 &  0.02 \\
            \PegaseHuge &  0.02 &  0.00 &  0.00 &   0.00 \\
            \GOC &  0.04 &  0.00 &  0.00 &  0.00 \\
            \MSR &  0.00 &   0.00 &   0.00 &   0.00 \\
            \LYON &   0.01 &   0.00 &   0.00 &   0.00 \\
            \FRANCE &  0.01 &   0.00 &  0.00 &  0.00 \\
            \bottomrule
        \end{tabular}
        \caption{Active Power Generation: MAE before and after running load flow for OPF.}
    \label{tab:OPF_load_flow_mae_before_after}
\end{table}

\begin{table}[!t]
        \centering
        \begin{tabular}{lrrrr}
            \toprule
            \textbf{MAE(pu)} &   DC &   LLOA$_1$ &  LLQCP &   LLLF \\
            \midrule
            \PegaseSmall &  0.06 &   0.00 &   0.00 &   0.00 \\
            \PegaseMedium &  0.05 &   0.00 &   0.00 &  0.01 \\
            \SDET &   0.02 &   0.00 &   0.00 & 0.00\\
            \RTE &  0.12 &  0.09 & 0.09 &  0.09 \\
            \PegaseLarge &  0.07 &  0.03 &   0.02 &  0.02 \\
            \MSR &  0.00 &   0.00 &   0.00 &   0.00 \\
            \LYON &   0.01 &   0.00 &   0.00 &   0.00 \\
            \FRANCE &  0.01 &   0.00 &  0.00 &  0.00 \\
            \bottomrule
        \end{tabular}
        \caption{Active Power Generation: MAE before and after running load flow for SCED.}
    \label{tab:SCED_load_flow_mae_before_after}
\end{table}

\begin{table}[!t]
        \centering
        \begin{tabular}{lrrrr}
            \toprule
            \textbf{Obj Gap(\%)} &   DC &   LLOA$_1$ &  LLQCP &   LLLF \\
            \midrule
            \PegaseSmall &  0.13 &   0.05 &   0.02 &   0.05 \\
            \PegaseMedium &  0.18 &   0.20 &   0.02 &  0.17 \\
            \SDET &   0.06 &   0.04 &   0.03 & 0.04\\
            \RTE &  0.16 &  0.11 & 0.11 &  0.11 \\
            \PegaseLarge &  0.17 &  0.16 &   0.02 &  0.15 \\
            \PegaseHuge &  0.06 &  0.04 &  0.01 &   0.04 \\
            \GOC &  0.01 &  0.00 &  0.00 &  0.00 \\
            \MSR &  0.01 &   0.00 &   0.00 &   0.00 \\
            \LYON &   0.05 &   0.01 &   0.01 &   0.01 \\
            \FRANCE &  0.02 &   0.01 &  0.00 &  0.01 \\
            \bottomrule
        \end{tabular}
        \caption{Active Power Generation: MAE with ac after running load flow for OPF.}
    \label{tab:OPF_load_flow_mae_compared_with ac}
\end{table}

\begin{table}[!t]
        \centering
        \begin{tabular}{lrrrr}
            \toprule
            \textbf{MAE(pu)} &   DC &   LLOA$_1$ &  LLQCP &   LLLF \\
            \midrule
            \PegaseSmall &  0.12 &   0.03 &   0.03 &   0.02 \\
            \PegaseMedium &  0.13 &   0.07 &   0.02 &  0.07 \\
            \SDET &   0.04 &   0.03 &   0.03 & 0.03\\
            \RTE &  0.16 &  0.12 & 0.11 &  0.12 \\
            \PegaseLarge &  0.18 &  0.13 &   0.05 &  0.13 \\
            \MSR &  0.02 &   0.00 &   0.00 &   0.00 \\
            \LYON &   0.03 &   0.02 &   0.02 &   0.03 \\
            \FRANCE &  0.02 &   0.01 &  0.00 &  0.01 \\
            \bottomrule
        \end{tabular}
        \caption{Active Power Generation: MAE with ac after running load flow for SCED.}
    \label{tab:SCED_load_flow_mae_compared_with ac}
\end{table}